\documentclass[11pt,a4paper]{article}
\usepackage{mathrsfs,amsfonts}
\usepackage[centertags]{amsmath}
\usepackage{amsthm,latexsym,graphicx,paralist}
\usepackage[all]{xy}
\newtheorem{theorem}{Theorem}[section]
\newtheorem{proposition}[theorem]{Proposition}
\newtheorem{lemma}[theorem]{Lemma}
\newtheorem{corollary}[theorem]{Corollary}

\newtheorem{definition}[theorem]{Definition}
\def\qed{\ifvmode\removelastskip\fi
{\unskip\nobreak\hfil\penalty50\hbox{}\nobreak\hfil \hbox{\vrule
height1.2ex width1.2ex}\parfillskip=0pt \finalhyphendemerits=0
\par \smallskip}}
\newcommand\rmd{\textup{d}}

\title{Geometric interpretations of the symmetric product in affine
differential geometry}
\author{\textsc{Mar\'ia Barbero-Li\~n\'an}\thanks{Work performed while a postdoctoral fellow at
Queen's University.} \\ Institute for the Mathematical Sciences
(CSIC-UAM-UC3M-UCM)\\ C/Nicol\'as Cabrera 13-15, \\ 28049 Madrid,
Spain \\ Email: mbarbero@icmat.es \and \textsc{Andrew D.\ Lewis}
\\ Department of Mathematics and Statistics \\ Queen's University\\
Kingston, ON K7L 3N6, Canada. \\ Email: andrew@mast.queensu.ca \\
}
\begin{document}
\maketitle

\begin{abstract}
The symmetric product of vector fields on a manifold arises when one studies
the controllability of certain classes of mechanical control systems.  A
geometric description of the symmetric product is provided using parallel
transport, along the lines of the flow interpretation of the Lie bracket.
This geometric interpretation of the symmetric product is used to provide an
intrinsic proof of the fact that the distributions closed under the symmetric
product are exactly those distributions invariant under the geodesic flow.
\end{abstract}

\vspace{4mm}


\section{Introduction}

Given an affine connection $\nabla$ on a manifold $M$\@, the corresponding
\textbf{symmetric product} is simply given by
\begin{equation*}
\langle X:Y\rangle=\nabla_XY+\nabla_YX.
\end{equation*}
The symmetric product for Levi-Civita connections appeared for first time
in~\cite{PEC:81} in the study of gradient systems.  This product appeared
again in~\cite{ADL/RMM:95c} where it was used to characterize the
controllability of a large class of mechanical control systems.  Since then,
the symmetric product has been widely used to solve control theoretic
problems for mechanical systems, such as motion
planning~\cite{FB/NEL/ADL:97,MK/JEM:10}\@,
trackability~\cite{MBL/MS:10,FB/ADL:04}\@, and so on.  We refer
to~\cite{FB/ADL:04} as a general reference for control theory for mechanical
systems.

The symmetric product has an interesting interpretation similar to that for
the Lie bracket as it relates to integrable distributions.  Let us recall the
result from~\cite{ADL:96b,ADL:96a}\@.  We say that a distribution
$\mathcal{D}$ on $M$ is \textbf{geodesically invariant} under an affine
connection $\nabla$ on $M$ if, as a submanifold of $TM$\@, $\mathcal{D}$ is
invariant under the geodesic spray associated with $\nabla$\@.  One can then
show that a distribution is geodesically invariant if and only if the
symmetric product of any $\mathcal{D}$-valued vector fields is again a
$\mathcal{D}$-valued vector field.  We provide an intrinsic proof of this
result in Section~\ref{sec:geodesic-invariance}\@.

Now, for the Lie bracket, one has the well-known formula
\begin{equation}\label{eq:Lie-bracket}
[X,Y](x)=\frac{1}{2}\left.\frac{\rmd^2}{\rmd t^2}\right|_{t=0}
\Phi^Y_{-t}\circ\Phi^X_{-t}\circ\Phi^Y_t\circ\Phi^X_t(x),
\end{equation}
where $\Phi^X_t$ denotes the flow of
$X$~\cite[Proposition~4.2.34]{RA/JEM/TSR:88}\@.  In this paper we provide for
the first time a similar formula for the symmetric product, using parallel
transport.  This is a novel interpretation.  Moreover, we use our
interpretation of the symmetric product to provide a coordinate-free proof of
the theorem on geodesic invariance mentioned in the preceding paragraph.  The
original proof in \cite{ADL:96b,ADL:96a} uses coordinates, and we refer
to~\cite{AB:07,AB:10} for an intrinsic proof using the bundle of linear
frames.

Let us provide an outline of the paper.  In Section~\ref{sec:background} we
provide our differential geometric notation and recall some facts that we
shall use in the paper.  One of the features of the paper is that it makes
essential and novel use of the
Baker\textendash{}Campbell\textendash{}Hausdorff formula and we review this
in Section~\ref{sec:BCH}\@.  In Section~\ref{sec:SymProd} we give various
infinitesimal descriptions of the symmetric product, see
Theorem~\ref{the:infinitesimal-symprod}\@.  In
Section~\ref{sec:geodesic-invariance} we use our infinitesimal descriptions
of the symmetric product to prove the geodesic invariance
theorem~\cite{ADL:96b,ADL:96a} mentioned above.  One of the contributions of
the paper is to give only intrinsic, coordinate-free characterizations and
proofs, and as a result there are many calculations in the paper that may be
of independent interest.  In particular, as mentioned above, we make use of
the Baker\textendash{}Campbell\textendash{}Hausdorff formula in a novel way
in a few places.

\section{Notation, background, and preliminary
constructions}\label{sec:background}

In this section we recall the basic facts about affine connections and
tangent bundles that will be important for us.  Some of our constructions are
presented in detail since we give\textemdash{}for the first time as far as we
are aware\textemdash{}some intrinsic definitions and proofs that are
well-known using coordinates.

Here is the notation we shall use in the paper.  By $\textup{Id}_S$ we denote
the identity map of a set $S$\@.  By $\mathbb{Z}_{\ge0}$ and $\mathbb{R}$ we
denote the set of nonnegative integers and real numbers, respectively.  For
the most part, we shall adopt the differential geometric conventions
of~\cite{RA/JEM/TSR:88}\@.  We shall assume all manifolds are paracompact,
Hausdorff, and of class $\mathcal{C}^\infty$\@.  All maps and geometric
objects will be assumed to be of class $\mathcal{C}^\infty$\@, and we shall
frequently use the word ``smooth'' to mean of class $\mathcal{C}^\infty$\@.
The set of smooth functions on a manifold $M$ is denoted by
$\mathcal{C}^\infty(M)$\@.  For a manifold $M$\@, its tangent bundle will be
denoted by $\tau_M\colon TM\rightarrow M$\@.  If $f\colon M\rightarrow N$ is
a map, its derivative is denoted by $Tf\colon TM\rightarrow TN$\@, and $T_xf$
denotes the restriction of $f$ to the tangent space $T_xM$\@.  The flow of a
vector field $X$ is denoted by $\Phi^X_t$\@,~i.e.,~the integral curve of $X$
through $x$ is $t\mapsto\Phi^X_t(x)$\@.  We shall suppose that all vector
fields are complete, and leave to the reader the task of modifying proofs to
account for the case where flows are defined on subintervals of
$\mathbb{R}$\@.  If $\pi\colon E\rightarrow M$ is a vector bundle over $M$\@,
we denote by $\Gamma^\infty(E)$ the set of smooth sections of $E$\@.
Sometimes it will be convenient to denote the zero vector in the fiber $E_x$
as $0_x$\@.  If $X\in\Gamma^\infty(TM)$ is a vector field and if $\Phi\colon
M\to M$ is a diffeomorphism, the pull-back of $X$ by $\Phi$ is given by
\begin{equation*}
\Phi^*X=T\Phi^{-1}\circ X\circ\Phi.
\end{equation*}
For a vector field $X\in\Gamma^\infty(TM)$ and for a function
$f\in\mathcal{C}^\infty(M)$\@, we denote by $\mathcal{L}_Xf$ the Lie
derivative of $f$ with respect to $f$\@.

\subsection{The Baker\textendash{}Campbell\textendash{}Hausdorff
formula}\label{sec:BCH}

One of the features of our presentation is that we use the
Baker\textendash{}Campbell\textendash{}Hausdorff (BCH) formula, as enunciated
in~\cite{RSS:87}\@, to evaluate compositions of flows in a crucial way in a
few places.  In this section we quickly review this formula.

The BCH formula provides a formula for the ``product of exponentials'' in a
Lie algebra in terms of brackets of the quantities being exponentiated.
First we recall the formal version of the formula, following~\cite{JPS:92}\@.
Consider a finite set $\boldsymbol{\xi}=\{\xi_1,\dots,\xi_p\}$ of
indeterminates and let $\hat{A}(\boldsymbol{\xi})$ be the
$\mathbb{R}$-algebra of formal power series in these indeterminates.  To be
clear about this, let $V(\boldsymbol{\xi})$ be the free $\mathbb{R}$-vector
space generated by $\boldsymbol{\xi}$\@.  Thus an element $\zeta\in
V(\boldsymbol{\xi})$ is a map $\zeta\colon\boldsymbol{\xi}\to\mathbb{R}$\@,
and the set of such maps is equipped with the pointwise operations of
addition and scalar multiplication.  For $k\in\mathbb{Z}_{\ge0}$\@, let
$T^k(V(\boldsymbol{\xi}))$ be the $k$th tensor power of
$V(\boldsymbol{\xi})$\@.  Then $\hat{A}(\boldsymbol{\xi})=
\prod_{k\in\mathbb{Z}_{\ge0}}T^k(V(\boldsymbol{\xi}))$ is the direct product.
Thus an element of $\hat{A}(\boldsymbol{\xi})$ is a map
\begin{equation*}
\alpha\colon\mathbb{Z}_{\ge0}\to
\cup_{k\in\mathbb{Z}_{\ge0}}T^k(V(\boldsymbol{\xi}))
\end{equation*}
such that $\alpha(k)\in T^k(V(\boldsymbol{\xi}))$\@.  The $\mathbb{R}$-vector
space $\hat{A}(\boldsymbol{\xi})$ is an algebra with the tensor product as
the product.  This algebra then has the natural Lie algebra structure given
by commutation: $[\alpha,\beta]=\alpha\beta-\beta\alpha$\@.  By
$\hat{L}(\boldsymbol{\xi})$ we denote the Lie subalgebra of
$\hat{A}(\boldsymbol{\xi})$ generated by the indeterminates
$\{\xi_1,\dots,\xi_k\}$\@.  Thus, formally, elements of
$\hat{L}(\boldsymbol{\xi})$ are $\mathbb{R}$-linear combinations of Lie
brackets of the indeterminates.  Let $L(\boldsymbol{\xi})$ be the Lie
subalgebra of $\hat{L}(\boldsymbol{\xi})$ having components in only finitely
many $T^k(V(\boldsymbol{\xi}))$\@,~i.e.,~the free Lie algebra generated by
the indeterminates $\boldsymbol{\xi}$\@.  One can then define a map
$\exp\colon\hat{L}(\boldsymbol{\xi})\to\hat{A}(\boldsymbol{\xi})$ by the
usual formal series expression:
\begin{equation*}
\exp(\alpha)=\sum_{k=0}^\infty\frac{\alpha^k}{k!}.
\end{equation*}
The formal Baker\textendash{}Campbell\textendash{}Hausdorff formula is then
the unique map
\begin{equation*}
\textup{BCH}\colon\underbrace{\hat{L}(\boldsymbol{\xi})\times\dots\times
\hat{L}(\boldsymbol{\xi})}_{k\ \textrm{copies}}\to\hat{L}(\boldsymbol{\xi})
\end{equation*}
satisfying
\begin{equation*}
\exp(\alpha_1)\cdots\exp(\alpha_k)=
\exp(\textup{BCH}(\alpha_1,\dots,\alpha_k)).
\end{equation*}
The component of $\textup{BCH}(\alpha_1,\dots,\alpha_k)$ in
$T^m(V(\boldsymbol{\xi}))$ we denote by
$\textup{BCH}_m(\alpha_1,\dots,\alpha_k)$\@, and we note that
\begin{equation}\label{eq:BCH12}
\begin{aligned}
\textup{BCH}_1(\alpha_1,\dots,\alpha_k)=&\;\alpha_1+\dots+\alpha_k,\\
\textup{BCH}_2(\alpha_1,\dots,\alpha_k)=&\;
\frac{1}{2}\sum_{\substack{a,b\in\{1,\dots,k\}\\a<b}}[\alpha_a,\alpha_b].
\end{aligned}
\end{equation}

Now let us recall what can be said about the BCH formula where the
indeterminates are vector fields $X_1,\dots,X_k$ on a manifold $M$\@.
The vector fields $X_1,\dots,X_k$ define a map
$\phi\colon\{\xi_1,\dots,\xi_k\}\to\Gamma^\infty(TM)$ by $\phi(\xi_j)=X_j$\@,
$j\in\{1,\dots,k\}$\@.  Since $\Gamma^\infty(TM)$ is a Lie algebra, there
exists a unique extension, which we also denote by $\phi$\@, from
$L(\boldsymbol{\xi})$ to $\Gamma^\infty(TM)$\@.  It is not generally the case
that the infinite formal series defining $\textup{BCH}$ converges (it does in
the real analytic case) in any reasonable topology on $\Gamma^\infty(TM)$\@,
but in~\cite{RSS:87} there are useful asymptotic formulae.  For our purposes,
these amount to the following.  For each $m\in\mathbb{Z}_{\ge0}$
\begin{equation}\label{eq:asympBCH}
\Phi^{X_k}_{t_k}\cdots\Phi^{X_1}_{t_1}(x)=
\sum_{j=1}^m\Phi^{\phi(\textup{BCH}_j(t_1X_1,\dots,t_kX_k))}_1(x)+
O((|t_1|+\dots+|t_k|)^{m+1})
\end{equation}
(here and subsequently, for brevity we denote composition of flows with
juxtaposition).  It is this formula that we shall use below.

\subsection{Tangent bundle geometry}\label{subsec:tangent}

In this section we review some well-known constructions concerning tangent
bundles.

We recall the definition of the vertical lift, which we regard as a vector
bundle map $\textup{vlft}\colon TM\oplus_MTM\to TTM$ as follows.  Let $x\in
M$ and let $v_x,w_x\in T_xM$\@.  The \textbf{vertical lift} of $u_x$ to $v_x$
is given by
\begin{equation}\label{eq:Vdef}
\textup{vlft}(v_x,u_x)=\left.\frac{\rmd}{\rmd t}\right|_{t=0}(v_x+tu_x).
\end{equation}
One easily verifies that the following diagram commutes:
\begin{equation*}\label{eq:Vdiagram}
\xymatrix{{TM\oplus_MTM}\ar[rr]^{\textup{vlft}}\ar[rd]_{\tau_M\oplus\tau_M}&&
{TTM}\ar[ld]^{T\tau_M\circ\tau_M}\\&{M}}
\end{equation*}
The image of $\textup{vlft}$ is the \textbf{vertical subbundle}
$VTM=\ker(T\tau_M)$\@.  For a vector field $X$ on $M$\@, the \textbf{vertical
lift} of $X$ is the vector field $X^V$ on $TM$ given by
$X^V(v_x)=\textup{vlft}(v_x,X(x))$\@.  It is evident that
\begin{equation}\label{eq:vlft-flow}
\Phi^{X^V}_t(v_x)=v_x+tX(x).
\end{equation}
In Section~\ref{subsec:dtb} we shall see that the double tangent bundle $TTM$
has two natural vector bundle structures, one for the vector bundle
$\tau_{TM}\colon TTM\to TM$ (called the primary vector bundle with the vector
bundle operations denoted with a subscript ``$1$'') and one for $T\pi\colon
TTM\to TM$ (called the secondary vector bundle with the vector bundle
operations denoted with a subscript ``$2$'').  The vertical lift interacts
with these two vector bundle structures differently in each component.
Indeed, the following diagrams commute
\begin{equation*}
\xymatrix{{TM\oplus_MTM}\ar[rr]^{\textup{vlft}}\ar[dr]_{\textup{pr}_1}&&
{TTM}\ar[dl]^{\tau_{TM}}\\&{TM}}\qquad
\xymatrix{{TM\oplus_MTM}\ar[rr]^{\textup{vlft}}\ar[dr]_{\textup{pr}_2}&&
{TTM}\ar[dl]^{T\tau_M}\\&{TM}}
\end{equation*}
This means that
\begin{equation}\label{eq:Vdiagrams}
\begin{gathered}
\textup{vlft}(v_1+v_2,u)=\textup{vlft}(v_1,u)+_2\textup{vlft}(v_2,u),\quad
\textup{vlft}(av,u)=a\cdot_2\textup{vlft}(v,u)\\
\textup{vlft}(v,u_1+u_2)=\textup{vlft}(v,u_1)+_1\textup{vlft}(v,u_2),\quad
\textup{vlft}(v,au)=a\cdot_1\textup{vlft}(v,u).
\end{gathered}
\end{equation}

We recall that, given $X\in\Gamma^\infty(TM)$\@, the \textbf{complete lift}
of $X$ is the vector field $X^C\in\Gamma^\infty(TTM)$ defined by
\begin{equation*}
X^C(v_x)=\left.\frac{\rmd}{\rmd t}\right|_{t=0}T_x\Phi^X_t(v_x).
\end{equation*}
Evidently,
\begin{equation}\label{eq:XC-flow}
\Phi^{X^C}_t(v_x)=T_x\Phi^X_t(v_x).
\end{equation}
Let us determine another useful characterisation of the flow of the complete
lift.
\begin{lemma}\label{lem:XC-flow}
Let\/ $X\in\Gamma^\infty(TM)$\@.  Let\/ $x_0\in M$ and\/ $v_0\in T_{x_0}M$\@.
Let\/ $J\subset\mathbb{R}$ be an interval for which\/ $0\in\textup{int}(J)$
and let\/ $\gamma\colon J\to M$ be a differentiable curve such that\/
$\gamma'(0)=v_0$\@.  Let\/ $I\subset\mathbb{R}$ be an interval.  Define\/
$\sigma(s,t)=\Phi^X_t(\gamma(s))$ for\/ $(s,t)\in J\times I$ and define a
vector field\/ $V_\sigma(t)=\left.\frac{\rmd}{\rmd
s}\right|_{s=0}\sigma(s,t)$ along the integral curve of\/ $X$ through\/
$x_0$\@.  Then the integral curve of\/ $X^C$ through\/ $v_0$ is\/ $t\mapsto
V_\sigma(t)$\@.
\begin{proof}
This is a simple computation:
\begin{equation*}
V_\sigma(t)=\left.\frac{\rmd}{\rmd s}\right|_{s=0}\Phi^X_t(\gamma(s))=
T_{\gamma(0)}\Phi^X_t(\gamma'(0))=T_{x_0}\Phi^X_t(v_0)=
\Phi^{X^C}_t(v_0)
\end{equation*}
using~\eqref{eq:XC-flow}\@.
\end{proof}
\end{lemma}

A consequence of the lemma is that the flow of $X^C$ is that of a linear
vector field, and so, by definition of a linear vector
field~\cite[\S47.9]{IK/PWM/JS:93}\@, $X^C$ is a vector bundle morphism
according to the following diagram:
\begin{equation}\label{eq:XCdiagram}
\xymatrix{{TM}\ar[r]^{X^C}\ar[d]_{\tau_M}&{TTM}\ar[d]^{T\pi}\\
{M}\ar[r]_{X}&{TM}}
\end{equation}

\subsection{The double tangent bundle}\label{subsec:dtb}

In this section we review some of the structure of the double tangent bundle
of a manifold.  We shall make great use of some of the constructions in this
section in our intrinsic constructions to follow.  Parts of the intrinsic
treatment we give of the canonical tangent bundle involution are, as far as
we know, new.

We begin by recalling the two vector bundle structures for $TTM$\@, as we
shall use both.  The double tangent bundle is represented naturally as a
vector bundle over $\tau_M\colon TM\to M$ in the following two ways:
\begin{equation}\label{eq:TTMdouble}
\xymatrix{{TTM}\ar[r]^{\tau_{TM}}\ar[d]_{T\tau_M}&{TM}\ar[d]^{\tau_M}\\
{TM}\ar[r]_{\tau_M}&{M}}\qquad\qquad
\xymatrix{{TTM}\ar[r]^{T\tau_M}\ar[d]_{\tau_{TM}}&{TM}\ar[d]^{\tau_M}\\
{TM}\ar[r]_{\tau_M}&{M}}
\end{equation}
The vector bundle on the left we call the \textbf{primary vector bundle} and
that on the right we call the \textbf{secondary vector bundle}\@.  We shall
need to introduce notation for the different vector bundle operations.  If
$u,v\in TTM$ satisfy $\tau_{TM}(u)=\tau_{TM}(v)$\@, then the sum of $u$ and
$v$ and the scalar multiple of $u$ by $a\in\mathbb{R}$ in the primary vector
bundle are denoted by $u+_1v$ and $a\cdot_1u$\@, respectively.  If $u,v\in
TTM$ satisfy $T\tau_M(u)=T\tau_M(v)$\@, then the sum of $u$ and $v$ and the
scalar multiple of $u$ by $a\in\mathbb{R}$ in the secondary vector bundle are
denoted by $u+_2v$ and $a\cdot_2u$\@, respectively.  For the vector bundle
$\tau_{TM}\colon TTM\to TM$\@, the vector bundle structure is the usual
tangent bundle structure.  We describe the vector bundle structure for
$T\tau_M\colon TTM\to TM$ as follows.  First note that the diagram
\begin{equation*}
\xymatrix{{TM}\ar[r]^{TX}\ar[d]_{\tau_M}&{TTM}\ar[d]^{\tau_{TM}}\\
{M}\ar[r]_{X}&{TM}}
\end{equation*}
commutes for a vector field $X$\@, giving $TX$ as a vector bundle mapping
over $X$\@.  Thus the map $X\mapsto TX$ is a morphism of the secondary vector
bundle structure.  Now let $u,v\in TTM$ be such that $w\doteq
T\tau_M(u)=T\tau_M(v)$\@.  We consider two cases.
\begin{enumerate}
\item $w\not=0$\@: Let $U,V\in\Gamma^\infty(TM)$ be such that $TU(w)=u$ and
$TV(w)=v$\@.  We then have
\begin{equation}\label{eq:secondary}
u+_2v=T(U+V)(w),\quad a\cdot_2u=T(aU)(w).
\end{equation}
\item $w=0$\@: In this case $u$ and $v$ are vertical.  So we let
$U,V\in\Gamma^\infty(TM)$ be such that $u=U^V\circ\tau_{TM}(u)$ and
$v=V^V\circ\tau_{TM}(v)$\@.  We then have
\begin{equation}\label{eq:secondarya}
u+_2v=(U+V)^V(\tau_{TM}(u)+\tau_{TM}(v)),\quad
a\cdot_2u=(aU)^V(a\tau_{TM}(u)).
\end{equation}
\end{enumerate}
We can say, motivated by~\eqref{eq:secondary}\@, that the secondary vector
bundle structure is the derivative of the vector bundle structure for
$\tau_M\colon TM\to M$\@.

The diagrams~\eqref{eq:TTMdouble} give a \textbf{double vector bundle} as
introduced in~\cite{JP:74}\@, and studied subsequently by many authors;
see~\cite[Chapter~9]{KCHM:05} for a general reference.  A consequence of this
structure is the following result that captures how the two vector bundle
structures are related.
\begin{lemma}\label{lem:intertwine}
Let\/ $u,v,w,z\in TTM$ satisfy
\begin{equation*}
T\tau_M(u)=T\tau_M(v),\quad T\tau_M(w)=T\tau_M(z),\quad
\tau_{TM}(u)=\tau_{TM}(w),\quad\tau_{TM}(v)=\tau_{TM}(z)
\end{equation*}
and let\/ $a,b\in\mathbb{R}$\@.  Then the following statements hold:
\begin{compactenum}
\item $(u+_2v)+_1(w+_2z)=(u+_1w)+_2(v+_1z)$\@;
\item $a\cdot_1(u+_2v)=(a\cdot_1u)+_2(a\cdot_1v)$\@;
\item $a\cdot_2(u+_1w)=(a\cdot_2u)+_1(a\cdot_2w)$\@;
\item $a\cdot_1(b\cdot_2w)=b\cdot_2(a\cdot_1w)$\@.
\end{compactenum}
\end{lemma}

To understand how the two vector bundle structures for $TTM$ are related, we
shall use a particular representation of points in $TTM$\@.  Let $\rho$ be a
smooth map from a neighbourhood of $(0,0)\in\mathbb{R}^2$ to $M$.  We shall
use coordinates $(s,t)$ for $\mathbb{R}^2$\@.  For fixed $s$ and $t$ define
$\rho_s(t)=\rho^t(s)=\rho(s,t)$\@, We then denote
\begin{equation*}
\frac{\partial}{\partial t}\rho(s,t)=\frac{\rmd}{\rmd t}\rho_s(t)
\in T_{\rho(s,t)}M,\quad
\frac{\partial}{\partial s}\rho(s,t)=\frac{\rmd}{\rmd s}\rho^t(s)
\in T_{\rho(s,t)}M.
\end{equation*}
Note that
\begin{equation*}
s\mapsto\frac{\partial}{\partial t}\rho(s,t)
\end{equation*}
is a curve in $TM$ for fixed $t$\@.  The tangent vector field to this curve
we denote by
\begin{equation*}
s\mapsto\frac{\partial}{\partial s}\frac{\partial}{\partial t}\rho(s,t)
\in T_{\frac{\partial}{\partial t}\rho(s,t)}TM.
\end{equation*}
We belabour the development of the notation somewhat since these partial
derivatives are not the usual partial derivatives from calculus, although the
notation might make one think they are.  For example, we do not generally
have equality of mixed partials,~i.e.,~generally we have
\begin{equation*}
\frac{\partial}{\partial s}\frac{\partial}{\partial t}\rho(s,t)\not=
\frac{\partial}{\partial t}\frac{\partial}{\partial s}\rho(s,t).
\end{equation*}

Now let $\rho_1$ and $\rho_2$ be smooth maps from a neighbourhood of
$(0,0)\in\mathbb{R}^2$ to $M$\@.  We say two such maps are
\textbf{equivalent} if
\begin{gather*}
\frac{\partial}{\partial s}\frac{\partial}{\partial t}\rho_1(0,0)=
\frac{\partial}{\partial s}\frac{\partial}{\partial t}\rho_2(0,0).
\end{gather*}
To the equivalence classes of this equivalence relation, we associate points
in $TTM$ by
\begin{equation*}
[\rho]\mapsto\frac{\partial}{\partial s}\frac{\partial}{\partial t}\rho(0,0).
\end{equation*}
We easily verify that
\begin{equation}\label{eq:rho-project}
\tau_{TM}([\rho])=\frac{\partial}{\partial t}\rho(0,0),\quad
T\tau_M([\rho])=\frac{\partial}{\partial s}\rho(0,0).
\end{equation}

Next, using the preceding representation of points in $TTM$\@, we relate the
two vector bundle structures for $TTM$ by defining a canonical involution of
$TTM$\@.  This is a well-known object, of course.  Our development and use of
this involution differs a little from what one usually sees in that it is
entirely free from local coordinates.  If $\rho$ is a smooth map from a
neighbourhood of $(0,0)\in\mathbb{R}^2$ into $M$\@, define another such map
by $\bar{\rho}(s,t)=\rho(t,s)$\@.  We then define the \textbf{canonical
tangent bundle involution} as the map $I_M\colon TTM\to TTM$ defined by
$I_M([\rho])=[\bar{\rho}]$\@.  Clearly $I_M\circ I_M=\textup{Id}_{TTM}$\@.

An interesting and useful formula connecting the complete lift and the
canonical tangent bundle involution is the following.
\begin{lemma}\label{lem:XCIM}
For\/ $X\in\Gamma^\infty(TM)$\@, $X^C=I_M\circ TX$\@.
\begin{proof}
Let $v_x\in TM$ and let $\gamma$ be a curve for which $\gamma'(0)=v_x$\@.  As
in Lemma~\ref{lem:XC-flow}\@, define $\sigma(s,t)=\Phi^X_t(\gamma(s))$ so
that
\begin{equation*}
X^C(v_x)=\frac{\partial}{\partial s}\frac{\partial}{\partial t}\sigma(0,0).
\end{equation*}
Then $\bar{\sigma}(s,t)=\Phi^X_s(\gamma(t))$ and so
\begin{equation*}
I_M(X^C(v_x))=\frac{\partial}{\partial t}
\frac{\partial}{\partial s}\bar{\sigma}(0,0)=
\left.\frac{\rmd}{\rmd t}\right|_{t=0}X(\gamma(t))=T_xX(v_x),
\end{equation*}
as desired.
\end{proof}
\end{lemma}

We have seen in~\eqref{eq:secondary} above that the secondary vector bundle
structure can be defined using the tangent functor.  Referring
to~\eqref{eq:XCdiagram} we see that $X\mapsto TX$ is a morphism with respect
to the primary vector bundle structure.  By the preceding lemma, this gives
us a way of representing the primary vector bundle operations in $TTM$\@.
Indeed, if $u,v\in TTM$ satisfy $\tau_{TM}(u)=\tau_{TM}(v)\doteq w$\@, we
consider the following two cases.
\begin{enumerate}
\item $w\not=0$\@: In this case, via~\eqref{eq:secondary} and the preceding
lemma, let $U,V\in\Gamma^\infty(TM)$ be such that $U^C(w)=u$ and
$V^C(w)=v$\@.  Then we have
\begin{equation}\label{eq:primary}
u+_1v=(U+V)^C(w),\quad a\cdot_1u=(aU)^C(w).
\end{equation}
\item $w=0$\@: In this case, $u,v\in T_{0_x}TM$ for a suitable $x$\@.  We
note that
\begin{equation*}
T_0TM\simeq T_xM\oplus T_xM,
\end{equation*}
cf.~\cite[Lemma~6.33]{FB/ADL:04}\@.  Thus there exists
$U,V\in\Gamma^\infty(TM)$ such that
\begin{equation*}
u=T\pi(u)\oplus U(x),\quad v=T\pi(v)\oplus V(x).
\end{equation*}
We then have
\begin{equation}\label{eq:primarya}
u+_1v=(T\pi(u)+T\pi(v))\oplus(U+V)(x),\quad
a\cdot_1u=(aT\pi(u))\oplus(aU)(x).
\end{equation}
\end{enumerate}

The following result will be helpful, and is more or less clear given the
preceding discussion.
\begin{lemma}\label{lem:IMbundle}
The map\/ $I_M$ is a vector bundle isomorphism:
\begin{equation*}
\xymatrix{{TTM}\ar[rr]^{I_M}\ar[rd]_{T\tau_M}&&{TTM}\ar[ld]^{\tau_{TM}}\\
&{TM}}
\end{equation*}
\begin{proof}
A proof in natural coordinates is elementary.  We shall give an intrinsic
proof.

It is clear from~\eqref{eq:rho-project} and the relations
\begin{equation*}
\frac{\partial}{\partial t}\bar{\rho}(0,0)=
\frac{\partial}{\partial s}\rho(0,0),\quad
\frac{\partial}{\partial s}\bar{\rho}(0,0)=
\frac{\partial}{\partial t}\rho(0,0)
\end{equation*}
that the diagram in the statement of the lemma commutes.  Moreover,
it is also clear that $I_M$ is a bijection.  It thus remains to show
that it is a vector bundle map.  Let $u,v\in TTM$ be such that
$T\tau_M(u)=T\tau_M(v)\doteq w$\@.  We then consider two cases.

$w\not=0$\@: Let $U,V\in\Gamma^\infty(TM)$ be such that $TU(w)=u$ and
$TV(w)=v$\@.  Then, using Lemma~\ref{lem:XCIM} and
equations~\eqref{eq:secondary} and~\eqref{eq:primary}\@,
\begin{align*}
I_M(u+_2v)=&\;I_M\circ T(U+V)(w)=(U+V)^C(w)\\
=&\;U^C(w)+_1V^C(w)=I_M\circ TU(w)+_1I_M\circ TV(w)\\
=&\;I_M(u)+_1I_M(v)
\end{align*}
and
\begin{align*}
I_M(a\cdot_2u)=&\;I_M\circ T(aU)(w)=(aU)^C(w)\\
=&\;a\cdot_1U^C(w)=a\cdot_1I_M\circ TU(w)=a\cdot_1I_M(u),
\end{align*}
as desired in this case.

$w=0$\@: Let $x=\tau_M\circ T\tau_M(u)=\tau_M\circ T\tau_M(v)$\@.
Choose $U,U',V,V'\in\Gamma^\infty(TM)$ such that
\begin{equation*}
u=U^V(U'(x)),\quad v=V^V(V'(x)).
\end{equation*}
For $s\in\mathbb{R}$ define $U_s,V_s\in\Gamma^\infty(TM)$ by
\begin{equation*}
U_s=U'+sU,\quad V_s=V'+sV.
\end{equation*}
Define
\begin{equation*}
\rho(s,t)=\Phi^{U_s}_t(x),\quad\sigma(s,t)=\Phi^{V_s}_t(x),
\end{equation*}
and note that
\begin{equation*}
\frac{\partial}{\partial t}\rho(s,0)=U_s(x),\quad
\frac{\partial}{\partial t}\sigma(s,0)=V_s(x)
\end{equation*}
and so $[\rho]=u$ and $[\sigma]=v$\@.  Now we use the
Baker\textendash{}Campbell\textendash{}Hausdorff formula to get
\begin{equation*}
\bar{\rho}(s,t)=\Phi^{U'+tU}_s(x)=\Phi^{sU'+stU}_1(x)
=\Phi^{sU}_t\circ\Phi^{U'}_s(x)+O((|s|+|t|)^2).
\end{equation*}
Therefore,
\begin{equation*}
\frac{\partial}{\partial t}\bar{\rho}(s,0)=sU(\Phi^{U'}_s(x))
\end{equation*}
and so
\begin{equation*}
\frac{\partial}{\partial s}\frac{\partial}{\partial t}\bar{\rho}(0,0)=
(0,U(x))+(U'(x),0)\in T_0TM\simeq T_xM\oplus T_xM.
\end{equation*}
Thus we have
\begin{equation*}
I_M(u)=T\tau_M(u)\oplus U(x),
\end{equation*}
with a similar formula holding for $v$\@, of course.  Therefore,
\begin{align*}
I_M(u+_2v)=&\;I_M(U^V(U'(x))+_2V^V(V'(x)))\\
=&\;I_M((U+V)^V(U'(x)+V'(x)))\\
=&\;(U'(x)+V'(x),U(x)+V(x))\\
=&\;I_M(u)+_1I_M(v),
\end{align*}
using~\eqref{eq:secondarya}\@,~\eqref{eq:primarya}\@, and the preceding
calculations.  Similarly,
\begin{equation*}
I_M(a\cdot_2u)=I_M((aU)^V(aU'(x)))=(aU'(x),aU(x))=a\cdot_1I_M(u),
\end{equation*}
as desired.
\end{proof}
\end{lemma}

We close this section with a few technical lemmata that we will subsequently
use in the paper.
\begin{lemma}\label{lem:IMvlft}
If\/ $w\in TTM$ satisfies\/ $\tau_{TM}(w)=v$ and\/ $T\tau_M(w)=u$ and
if\/ $z\in T_xM$\@, then
\begin{equation*}
w+_2I_M\circ\textup{vlft}(u,z)=w+_1\textup{vlft}(v,z).
\end{equation*}
\begin{proof}
Let $U,V,Z\in\Gamma^\infty(TM)$ be vector fields for which
\begin{equation*}
U(x)=u,\quad V(x)=v,\quad Z(x)=z,
\end{equation*}
We consider two cases.

$u\not=0$\@: In this case, write $w=TW(u)$ for some vector field
$W\in\Gamma^\infty(TM)$\@.  Then $W(x)=V(x)$\@.  We compute
\begin{align*}
w+_2I_M\circ\textup{vlft}(u,z)=&\;TW(u)+_2I_M(Z^V(U(x)))\\
=&\;TW(u)+_2(U(x)\oplus Z(x))\\
=&\;TW(u)+_1Z^V(V(x)),
\end{align*}
as desired in this case.

$u=0$\@: Here we write $w=W^V(V(x))$ for an appropriate vector field $W$ on
$M$\@.  Then
\begin{align*}
w+_2I_M\circ\textup{vlft}(u,z)=&\;W^V(V(x))+_2I_M(Z^V(U(x)))\\
=&\;W^V(V(x))+_2(U(x)\oplus Z(x))\\
=&\;W^V(V(x))+_1Z^V(V(x)),
\end{align*}
giving the lemma.
\end{proof}
\end{lemma}

The proof of the following lemma is a specialization of the proof of
Lemma~6.19 in~\cite{IK/PWM/JS:93}\@.
\begin{lemma}\label{lem:[X,Y]V}
For\/ $X,Y\in\Gamma^\infty(TM)$ we have
\begin{equation*}
TY(X(x))-_1I_M\circ TX(Y(x))=\textup{vlft}(Y(x),[X,Y](x)).
\end{equation*}
\begin{proof}
We use the formula
\begin{equation*}
[X,Y](x)=\left.\frac{\rmd}{\rmd t}\right|_{t=0}(\Phi^X_t)^*Y(x),
\end{equation*}
\cite[Theorem~4.2.19]{RA/JEM/TSR:88}\@.  Note that the curve
\begin{equation*}
t\mapsto(\Phi^X_t)^*Y(x)
\end{equation*}
is a curve in $T_xM$ passing through $Y(x)$ at $t=0$\@, and so its
derivative with respect to $t$ at $t=0$ is a vertical tangent vector
in $T_{Y(x)}TM$. Note that $V_{Y(x)}TM\simeq T_xM$\@.  We calculate
\begin{align*}
\left.\frac{\rmd}{\rmd t}\right|_{t=0}(\Phi^X_t)^*Y(x)=&\;
\left.\frac{\rmd}{\rmd t}\right|_{t=0}T\Phi^X_{-t}\circ Y\circ\Phi^X_t(x)\\
=&\;-X^C\circ Y(x)+_1TY\circ X(x)\\
=&\;TY\circ X(x)-_1I_M\circ TX\circ Y(x),
\end{align*}
using Lemma~\ref{lem:XCIM} and~\eqref{eq:XC-flow}\@.
\end{proof}
\end{lemma}

\subsection{Affine differential geometry}

This section will be a very rapid overview of the affine differential
geometry we shall use in this paper.  We refer to~\cite{SK/KN:63a} for
details.

A $\mathcal{C}^\infty$-\textbf{affine connection} on a manifold $M$ assigns
to each pair $(X,Y)\in\Gamma^\infty(TM)\times\Gamma^\infty(TM)$ a vector
field $\nabla_XY\in\Gamma^\infty(TM)$\@, and the assignment satisfies
\begin{compactenum}
\item the map $(X,Y)\mapsto\nabla_XY$ is $\mathbb{R}$-bilinear,
\item $\nabla_{fX}Y =f\nabla_XY$\@, and
\item $\nabla_X(fY)=f\nabla_XY +(\mathcal{L}_X f)Y$
\end{compactenum}
for each $X,Y\in\Gamma^\infty(TM)$ and $f\in\mathcal{C}^\infty(M)$\@.  The
vector field $\nabla_XY$ is called the \textbf{covariant derivative} of $Y$
with respect to $X$\@.

As the expression $\nabla_XY$ is tensorial in $X$\@, it only depends on the
value of $X$ at the point $x$\@.  Hence, if $v_x\in T_xM$\@, we can define
\begin{equation*}
\nabla_{v_x} Y(x)=\nabla_XY(x)\in T_xM,
\end{equation*}
where $X$ is any $\mathcal{C}^\infty$-vector field such that $X(x)=v_x$\@.

Given an affine connection $\nabla$\@, there exists a complementary subbundle
$HTM$ of the vertical subbundle $VTM=\ker(T\tau_M)$\@,~i.e.,~$TTM=HTM\oplus
VTM$\@.  This complementary subbundle is called the \textbf{horizontal
subbundle} and is constructed as follows~\cite{IK/PWM/JS:93}\@.  We shall
first define a map $\textup{hlft}\colon TM\oplus_MTM\to TTM$\@.  Let $x\in M$
and $u,v\in T_xM$\@.  Let $X\in\Gamma^\infty(TM)$ be such that $X(x)=v$ and
define
\begin{equation}\label{eq:Hdef}
\textup{hlft}(v,u)=TX(u)-_1\textup{vlft}(v,\nabla_uX),
\end{equation}
where $\textup{vlft}$ is the vertical lift map from~\eqref{eq:Vdef}\@.  One
can easily check that $\textup{hlft}$ is indeed a vector bundle map according
to both of the following commuting diagrams:
\begin{equation*}
\xymatrix{{TM\oplus_MTM}\ar[rr]^{\textup{hlft}}\ar[dr]_{\textup{pr}_1}&&
{TTM}\ar[dl]^{\tau_{TM}}\\&{TM}}\qquad
\xymatrix{{TM\oplus_MTM}\ar[rr]^{\textup{hlft}}\ar[dr]_{\textup{pr}_2}&&
{TTM}\ar[dl]^{T\tau_M}\\&{TM}}
\end{equation*}
Thus
\begin{equation}\label{eq:Hdiagrams}
\begin{gathered}
\textup{hlft}(v_1+v_2,u)=\textup{hlft}(v_1,u)+_2\textup{hlft}(v_2,u),\quad
\textup{hlft}(av,u)=a\cdot_2\textup{hlft}(v,u)\\
\textup{hlft}(v,u_1+u_2)=\textup{hlft}(v,u_1)+_1\textup{hlft}(v,u_2),\quad
\textup{hlft}(v,au)=a\cdot_1\textup{hlft}(v,u).
\end{gathered}
\end{equation}
The horizontal subbundle is defined by
\begin{equation*}
H_{v_x}TM=\{\textup{hlft}(v_x,u_x)|\enspace u_x\in T_xM\}.
\end{equation*}
At each $v_x\in TM$\@, the linear map $T_{v_x}\tau_M\colon
T_{v_x}TM\rightarrow T_xM$\@, restricted to the horizontal subspace
$H_{v_x}TM$\@, is an isomorphism.  The inverse of this isomorphism,
applied to $u_x\in T_xM$\@, is the \textbf{horizontal lift} of $u_x$
to $v_x\in T_xM$\@:
\begin{equation*}
\textup{hlft}(v_x,u_x)=(T_{v_x}\tau_M|H_{v_x}TM)^{-1}(u_x).
\end{equation*}
The \textbf{horizontal lift} of the vector field
$X\in\Gamma^\infty(TM)$ is the vector field
$X^H\in\Gamma^\infty(TTM)$ defined by
$X^H(v_x)=\textup{hlft}(v_x,X(x))$\@.

The \textbf{torsion tensor} is denoted by $T$\@:
\begin{equation*}
T(X,Y)=\nabla_XY-\nabla_YX-[X,Y].
\end{equation*}
The canonical tangent bundle involution also provides an interesting and
useful way of characterising torsion-free affine connections.  The following
result appears in~\cite{RJF/HTL:99}\@, but with a coordinate proof.  We
provide an intrinsic proof that is quite a lot simpler than the proof
in~\cite{RJF/HTL:99}\@.
\begin{lemma}
With the notation preceding,
\begin{equation*}
\textup{hlft}(v_x,u_x)-_1I_M\circ\textup{hlft}(u_x,v_x)=
\textup{vlft}(v_x,T(v_x,u_x))
\end{equation*}
for all\/ $u_x,v_x\in T_xM$ and all\/ $x\in M$\@.  As a consequence, the
following statements are equivalent:
\begin{compactenum}
\item \label{pl:torsion1} $\nabla$ is torsion-free;
\item \label{pl:torsion2}
$\textup{hlft}(v_x,u_x)=I_M\circ\textup{hlft}(u_x,v_x)$ for all\/
$u_x,v_x\in T_xM$ and\/ $x\in M$\@;
\item \label{pl:torsion3} $I_M$ leaves the horizontal subbundle\/ $HTM\subset
TTM$ invariant.
\end{compactenum}
\begin{proof}
The first assertion of the lemma follows from Lemmata~\ref{lem:IMvlft}
and~\ref{lem:[X,Y]V} as follows:
\begin{align*}
\textup{hlft}&(X(x),Y(x))-_1I_M\circ\textup{hlft}(Y(x),X(x))\\
=&\;(TX(Y(x))-_1\textup{vlft}(X(x),\nabla_YX(x)))
-_1I_M(TY(X(x))-_1\textup{vlft}(Y(x),\nabla_XY(x)))\\
=&\;TX(Y(x))-_1\textup{vlft}(X(x),\nabla_YX(x))
-_1(I_M\circ TY(X(x))-_2I_M\circ\textup{vlft}(Y(x),\nabla_XY(x)))\\
=&\;TX(Y(x))-_1\textup{vlft}(X(x),\nabla_YX(x))
-_1(I_M\circ TY(X(x))-_1\textup{vlft}(X(x),\nabla_XY(x)))\\
=&\;TX(Y(x))-_1\textup{vlft}(X(x),\nabla_YX(x))\\
&\;-_1(TX(Y(x))+_1\textup{vlft}(X(x),[X,Y](x))-_1
\textup{vlft}(X(x),\nabla_XY(x)))\\
=&\;\textup{vlft}(X(x),-[X,Y](x)+\nabla_XY(x)-\nabla_YX(x))\\
=&\;\textup{vlft}(X(x),T(X(x),Y(x))),
\end{align*}
for vector fields $X$ and $Y$\@.

\eqref{pl:torsion1}$\implies$\eqref{pl:torsion2} This follows immediately
from the first assertion of the lemma.

\eqref{pl:torsion2}$\implies$\eqref{pl:torsion3} This is obvious.

\eqref{pl:torsion3}$\implies$\eqref{pl:torsion1} Let $w\in HTM$ so that
$I_M(w)\in HTM$\@.  Write $w=\textup{hlft}(v_x,u_x)$ for some $x\in M$ and
$u_x,v_x\in T_xM$\@.  Since
\begin{equation*}
\tau_{TM}(w)=v_x,\quad T\tau_M(w)=u_x,
\end{equation*}
we have
\begin{equation*}
\tau_{TM}(I_M(w))=u_x,\quad T\tau_M(I_M(w))=v_x.
\end{equation*}
Since $I_M(w)$ is horizontal, we must have $I_M(w)=\textup{hlft}(u_x,v_x)$\@.
It then immediately follows that $T=0$ from the first part of the proof.
\end{proof}
\end{lemma}

Given an interval $I\subset\mathbb{R}$ and a curve $\gamma\colon I\rightarrow
M$\@, a \textbf{vector field along $\gamma$} is a smooth map that assigns to
every $t\in I$ an element of $T_{\gamma(t)}M$\@.  If $Y\colon I\to TM$ is a
vector field along $\gamma$\@, it makes sense to define a
$\mathcal{C}^\infty$-vector field along $\gamma$ by
\begin{equation*}
I\ni t\mapsto\nabla_{\gamma'(t)}\overline{Y}(\gamma(t))\in T_{\gamma(t)}M,
\end{equation*}
where $\overline{Y}$ is a vector field for which
$Y(t)=\overline{Y}(\gamma(t))$\@.  This construction can be shown to be
independent of the extension of $Y$ to $\overline{Y}$\@.  A vector field $Y$
along $\gamma$ is \textbf{parallel} if $\nabla_{\gamma'(t)}Y(t)=0$ for each
$t\in I$\@.

The equation $\nabla_{\gamma'(t)}Y(t)=0$ can be regarded as a differential
equation for the vector field $Y$ along $\gamma$\@.  If the initial value $v$
of the vector field at $t_0\in I$ is given, the differential equation has a
unique solution $Y(t)$ for $t$ sufficiently close to $t_0$\@.  The map
$\tau_\gamma^{(t,t_0)}\colon T_{\gamma(t_0)}M \rightarrow T_{\gamma(t)}M$
that sends $v\in T_{\gamma(t_0)}M$ to the unique vector $Y(t)\in
T_{\gamma(t)}M$ defined by the solution to the initial value problem
\begin{equation*}
\nabla_{\gamma'(t)}Y(t)=0,\quad Y(t_0)=v,
\end{equation*}
is called the \textbf{parallel transport} along $\gamma$\@.  Note
that $\tau_\gamma^{(t,t_0)}$ is an isomorphism.  We recall from the
discussion in~\cite[page~114]{SK/KN:63a} the following formula:
\begin{equation}\label{eq:hlft-flow}
\Phi^{X^H}_t(v_x)=\tau^{(t,0)}_{\gamma}(v_x),
\end{equation}
where $\gamma$ is the integral curve of the vector field $X$ for
which $\gamma(0)=x$\@. The covariant derivative of $Y$ along $X$ can
also be described as follows:
\begin{equation}\label{eq:covartau}
\nabla_XY(x)=\left.\frac{\rmd}{\rmd t}
\right|_{t=0}\tau^{(0,t)}_\gamma(Y(\gamma(t))),
\end{equation}
where $\gamma$ is the integral curve of $X$ satisfying
$\gamma(0)=x$\@.

A \textbf{geodesic} of an affine connection $\nabla$ on $M$ is a
curve $\gamma\colon I \rightarrow M$ satisfying
$\nabla_{\gamma'(t)}\gamma'(t)=0$\@.  A geodesic can also be
described as a curve whose tangent vector field is parallel along
itself.  The geodesic equations give rise to a second-order vector
field $Z\in\Gamma^\infty(TTM)$ having the property that the integral
curves of $Z$ projected to $M$ by the natural tangent bundle
projection $\tau_M$ are geodesics of $\nabla$\@.  This vector field
$Z$ is called \textbf{geodesic spray} for $\nabla$\@.  The geodesic
spray can be defined using horizontal lifts as follows:
\begin{equation}\label{eq:Zdefn}
Z(v_x)=\textup{hlft}(v_x,v_x).
\end{equation}

Note that while parallel transport uses ``all'' of the information about an
affine connection, the geodesics do not, as they depend only on the symmetric
part of the Christoffel symbols.  This observation is made precise as
follows.  If $\nabla$ is an affine connection on $M$\@, then there exists a
unique torsion-free affine connection, denoted by $\overline{\nabla}$\@,
whose geodesics are exactly those of $\nabla$\@.  Explicitly,
\begin{equation}\label{eq:AffConn-Torsion}
\overline{\nabla}_XY=\nabla_XY-\frac{1}{2}T(X,Y),
\end{equation}
cf.~Propositions~7.9 and~7.10 in Chapter~III in~\cite{SK/KN:63a}\@.  Here $T$
is the torsion of $\nabla$\@.  It is possible to relate the parallel
transport of a connection and its torsion-free connection.
\begin{lemma}\label{lem:PT-SameGeod-DiffT}
Let\/ $\nabla$ be an affine connection on\/ $M$ with torsion\/ $T$ and let\/
$\overline{\nabla}$ be the corresponding zero-torsion affine connection.
Let\/ $\gamma$ be a geodesic for both\/ $\nabla$ and\/ $\overline{\nabla}$
with the same initial condition.  If\/ $V\in T_{\gamma(0)}M$ then
\begin{equation}\label{eq:Diff-Paral-Trans}
\tau^{(t,0)}_\gamma(V)-\overline{\tau}^{(t,0)}_\gamma(V)
=\overline{\tau}_\gamma^{(t,0)}\left(-\frac{1}{2}
\int^t_0\overline{\tau}_\gamma^{(0,s)}\left(
T(\gamma'(s),\tau_\gamma^{(s,0)}(V))\right)\,\rmd s\right),
\end{equation}
where\/ $\tau^{t,0}_\gamma$ (resp.~$\overline{\tau}_\gamma^{t,0})$ is the\/
$\nabla$ (resp.~$\overline{\nabla})$ parallel transport along\/ $\gamma$
from\/ $T_{\gamma(0)}M$ to\/ $T_{\gamma(t)}M$\@.
\begin{proof}
Let us abbreviate
\begin{equation*}
A_V(t)=-\frac{1}{2}
\int^t_0\overline{\tau}_\gamma^{(0,s)}\left(
T(\gamma'(s),\tau_\gamma^{(s,0)}(V))\right)\,\rmd s
\end{equation*}
so that
\begin{align*}
\frac{\rmd}{\rmd t}A_V(t)=&\;\overline{\tau}_\gamma^{(0,t)}
\left(-\frac{1}{2}T(\gamma'(t),\tau_\gamma^{(t,0)}(V))\right)\\
=&\;\overline{\tau}_\gamma^{(0,t)}
\left(\overline{\nabla}_{\gamma'(t)}\tau_\gamma^{(t,0)}(V)-
\nabla_{\gamma'(t)}\tau_\gamma^{(t,0)}(V)\right)\\
=&\;\overline{\tau}_\gamma^{(0,t)}
\left(\overline{\nabla}_{\gamma'(t)}\tau_\gamma^{(t,0)}(V)\right),
\end{align*}
using~\eqref{eq:AffConn-Torsion}\@.  We also compute
\begin{align*}
\frac{\rmd}{\rmd t}\left(\overline{\tau}_\gamma^{(0,t)}\circ
\tau_\gamma^{(t,0)}(V)-V\right)=&\;
\overline{\nabla}_{\gamma'(t)}\left(\overline{\tau}_\gamma^{(0,t)}\circ
\tau_\gamma^{(t,0)}(V)\right)\\
=&\;\overline{\tau}_\gamma^{(0,t)}\overline{\nabla}_{\gamma'(t)}
\tau_\gamma^{(t,0)}(V).
\end{align*}
Thus we have
\begin{equation*}
\frac{\rmd}{\rmd t}A_V(t)=
\frac{\rmd}{\rmd t}\left(\overline{\tau}_\gamma^{(0,t)}\circ
\tau_\gamma^{(t,0)}(V)-V\right).
\end{equation*}
Since
\begin{equation*}
A_V(t)|_{t=0}=\left.\left(\overline{\tau}_\gamma^{(0,t)}\circ
\tau_\gamma^{(t,0)}(V)-V\right)\right|_{t=0},
\end{equation*}
it follows that
\begin{equation*}
A_V(t)=\overline{\tau}_\gamma^{(0,t)}\circ\tau_\gamma^{(t,0)}(V)-V.
\end{equation*}
Rearranging gives the result.
\end{proof}
\end{lemma}

\section{Infinitesimal descriptions of the symmetric
product}\label{sec:SymProd}

Now we are ready to geometrically describe the symmetric product for vector
fields.  We shall provide four equivalent infinitesimal descriptions of the
symmetric product (some of which are related in elementary ways).  To do
this, we make use of the BCH formula.

Let $\nabla$ be an affine connection.  The \textbf{symmetric product} for
$\nabla$ of two vector fields $X$ and $Y$ on $M$ is defined as follows:
\begin{equation*}
\langle X\colon Y\rangle=\nabla_XY+\nabla_YX.
\end{equation*}
Our infinitesimal descriptions of the symmetric product, like that
of~\eqref{eq:Lie-bracket} for the Lie bracket, involve
concatenations of flows of vector fields.  Before we state the
results, let us give the various constructions we use.  We let
$\nabla$ be an affine connection on $M$ with $\overline{\nabla}$ the
associated zero-torsion connection.  We let
$X_1,X_2\in\Gamma^\infty(TM)$ and let $v_x\in TM$\@.  By $X_1^H$ and
$X_2^H$ we denote the horizontal lifts with respect to $\nabla$ and
by $\overline{X}_1^H$ and $\overline{X}_2^H$ we denote the
horizontal lifts with respect to $\overline{\nabla}$\@.  By $\eta_1$
and $\eta_2$ we denote the integral curves of $X_1$ and $X_2$\@,
respectively, through $x$\@.  We let $\tau_\gamma^{(t,0)}$ and
$\overline{\tau}_\gamma^{(t,0)}$ denote the parallel transport with
respect to $\nabla$ and $\overline{\nabla}$\@, respectively, along a
curve $\gamma$\@.  Now define four curves $\Upsilon_1$\@,
$\Upsilon_2$\@, $\Upsilon_3$\@, and $\Upsilon_4$ in $TM$ as follows:
\begin{align*}
\Upsilon_1(t)=&\;\Phi^{X_2^V}_{-t}\Phi^{X_1^H}_{-t}\Phi^{X_2^V}_{t}
\Phi^{X_1^H}_{t}\Phi^{X_1^V}_{-t}\Phi^{X_2^H}_{-t}\Phi^{X_1^V}_{t}
\Phi^{X_2^H}_{t}(v_x),\\
\Upsilon_2(t)=&\;\Phi^{X_2^V}_{-t}\Phi^{\overline{X}_1^H}_{-t}
\Phi^{X_2^V}_{t}\Phi^{\overline{X}_1^H}_{t}\Phi^{X_1^V}_{-t}
\Phi^{\overline{X}_2^H}_{-t}\Phi^{X_1^V}_{t}
\Phi^{\overline{X}_2^H}_{t}(v_x),\\
\Upsilon_3(t)=&\;\Phi^{X_2^V}_{-t}\tau_{\eta_1}^{(0,t)}
\Phi^{X_2^V}_{t}\tau_{\eta_1}^{(t,0)}\Phi^{X_1^V}_{-t}
\tau_{\eta_2}^{(0,t)}\Phi^{X_1^V}_{t}\tau_{\eta_2}^{(t,0)}(v_x),\\
\Upsilon_4(t)=&\;\Phi^{X_2^V}_{-t}\overline{\tau}_{\eta_1}^{(0,t)}
\Phi^{X_2^V}_{t}\overline{\tau}_{\eta_1}^{(t,0)}\Phi^{X_1^V}_{-t}
\overline{\tau}_{\eta_2}^{(0,t)}\Phi^{X_1^V}_{t}
\overline{\tau}_{\eta_2}^{(t,0)}(v_x).
\end{align*}

Before we state the main result in this section, we have the following lemma
that we shall use in its proof.  This formula appears, for example,
in~\cite{MC:00}\@.
\begin{lemma}\label{lem:crampin}
For vector fields\/ $X,Y\in\Gamma^\infty(TM)$ and for an affine connection\/
$\nabla$ on\/ $M$\@, we have\/ $(\nabla_XY)^V=[X^H,Y^V]$\@.
\begin{proof}
We use~\eqref{eq:Lie-bracket}\@:
\begin{equation*}
[X^H,Y^V](v_x)=-\frac{1}{2}\left.\frac{\rmd^2}{\rmd t^2}\right|_{t=0}
\Phi^{X^H}_{-t}\Phi^{Y^V}_{-t}\Phi^{X^H}_t\Phi^{Y^V}_t(v_x).
\end{equation*}
By equations~\eqref{eq:vlft-flow} and~\eqref{eq:hlft-flow} and by
linearity of parallel transport we compute
\begin{equation*}
\Phi^{X^H}_{-t}\Phi^{Y^V}_{-t}\Phi^{X^H}_t\Phi^{Y^V}_t(v_x)=
v_x-t(\tau^{(0,t)}_\eta(Y(\eta(t)))-Y(x)),
\end{equation*}
where $\eta$ is the integral curve of $X$ through $x$\@.  Note that
this is a curve in $T_xM$ and so its derivatives will be vertical
tangent vectors. We then have
\begin{equation*}
\left.\frac{\rmd}{\rmd t}\right|_{t=0}
\Phi^{X^H}_{-t}\Phi^{Y^V}_{-t}\Phi^{X^H}_t\Phi^{Y^V}_t(v_x)=0
\end{equation*}
and by~\eqref{eq:covartau}
\begin{equation*}
\left.\frac{\rmd^2}{\rmd t^2}\right|_{t=0}
\Phi^{X^H}_{-t}\Phi^{Y^V}_{-t}\Phi^{X^H}_t\Phi^{Y^V}_t(v_x)=
-2\textup{vlft}(v_x,\nabla_XY(x)),
\end{equation*}
from which the lemma immediately follows.
\end{proof}
\end{lemma}

With the preceding notation, we state the following theorem.
\begin{theorem}\label{the:infinitesimal-symprod}
With the notation of the preceding paragraph, if\/
$\Upsilon\in\{\Upsilon_1,\Upsilon_2,\Upsilon_3,\Upsilon_4\}$\@, then\/
$\Upsilon'(0)=0$ and
\begin{equation*}
\frac{1}{2}\left.\frac{\rmd^2}{\rmd t^2}\right|_{t=0}\Upsilon(t)=
\langle X_1:X_2\rangle^V(v_x).
\end{equation*}
\begin{proof}
Let us apply the BCH formulae~\eqref{eq:BCH12} to the concatenation of flows
defining $\Upsilon_1$\@.  It is immediately clear that
\begin{equation*}
\phi(\textup{BCH}_1(tX_2^H,tX_1^V,-tX_2^H,-tX_1^V,tX_1^H,
tX_2^V,-tX_1^H,-tX_2^V)=0.
\end{equation*}
Some bookkeeping and the fact the flows of vertically lifted vector fields
obviously commute gives
\begin{multline*}
\phi(\textup{BCH}_2(tX_2^H,tX_1^V,-tX_2^H,-tX_1^V,tX_1^H,
tX_2^V,-tX_1^H,-tX_2^V)\\=t^2([X_1^H,X_2^V]+[X_2^H,X_1^V]).
\end{multline*}
An application of~\eqref{eq:asympBCH} and Lemma~\ref{lem:crampin} now gives
\begin{align*}
\frac{1}{2}\left.\frac{\rmd^2}{\rmd t^2}\right|_{t=0}\Upsilon_1(t)=&\;
\frac{1}{2}\left.\frac{\rmd^2}{\rmd t^2}\right|_{t=0}
\Phi^{t^2\langle X_1:X_2\rangle^V}_{1}(v_x)\\
=&\;\frac{1}{2}\left.\frac{\rmd^2}{\rmd t^2}\right|_{t=0}
\Phi^{\langle X_1:X_2\rangle^V}_{t^2}(v_x)=\langle X_1:X_2\rangle^V(v_x).
\end{align*}

The same argument as above gives
\begin{equation*}
\frac{1}{2}\left.\frac{\rmd^2}{\rmd t^2}\right|_{t=0}\Upsilon_2(t)=
[\overline{X}_1^H,X_2^V](v_x)+[\overline{X}_2^H,X_1^V](v_x).
\end{equation*}
By~\eqref{eq:AffConn-Torsion} we have
\begin{equation*}
\overline{X}^H(v_x)=X^H(v_x)+\frac{1}{2}\textup{vlft}(v_x,T(X(x),v_x)).
\end{equation*}
As a result, one directly computes
\begin{equation*}
[\overline{X}^H,Y^V]=[X^H,Y^V]-\frac{1}{2}T(X,Y)^V.
\end{equation*}
Skew-symmetry of the torsion then gives
\begin{equation*}
[\overline{X}_1^H,X_2^V](v_x)+[\overline{X}_2^H,X_1^V](v_x)=
[X_1^H,X_2^V](v_x)+[X_2^H,X_1^V](v_x),
\end{equation*}
and from this we arrive at
\begin{equation*}
\frac{1}{2}\left.\frac{\rmd^2}{\rmd t^2}\right|_{t=0}\Upsilon_2(t)=
\langle X_1:X_2\rangle^V(v_x).
\end{equation*}

Given this formula and the results from the first part of the proof,
we immediately have from~\eqref{eq:hlft-flow}
\begin{equation*}
\frac{1}{2}\left.\frac{\rmd^2}{\rmd t^2}\right|_{t=0}\Upsilon_3(t)=
\frac{1}{2}\left.\frac{\rmd^2}{\rmd t^2}\right|_{t=0}\Upsilon_4(t)=
\langle X_1:X_2\rangle^V(v_x).
\end{equation*}
as desired.
\end{proof}
\end{theorem}

The following corollary gives a geometric interpretation of what is going on
with the composition of flows in the preceding theorem.
\begin{corollary}\label{cor:simple-infinitesimal}
Let\/ $X_1,X_2\in\Gamma^\infty(TM)$\@, let\/ $\nabla$ be an affine connection
on\/ $M$\@, and let\/ $x\in M$\@.  If\/ $\Upsilon_1=\Upsilon_3$ are defined
as preceding Theorem~\ref{the:infinitesimal-symprod} while taking\/
$v_x=0$\@, we have
\begin{equation}\label{eq:symprod0x}
\Upsilon_1(t)=\Upsilon_3(t)=
v_x+t\left(\tau_{\eta_2}^{(0,t)}(X_1(\eta_2(t)))-X_1(x)+
\tau_{\eta_1}^{(0,t)}(X_2(\eta_1(t)))-X_2(x)\right).
\end{equation}
\begin{proof}
The idea is the same, but only a little longer to carry out, as the proof of
Lemma~\ref{lem:crampin}\@.
\end{proof}
\end{corollary}

The upshot of the corollary is that the conclusion of
Theorem~\ref{the:infinitesimal-symprod} can be rendered a little more
transparent since it is more or less obvious that the first derivative of the
right-hand side of~\eqref{eq:symprod0x} is zero and that the second
derivative is twice the symmetric product.

There exists another infinitesimal description of the symmetric product along
the same lines as that of Theorem~\ref{the:infinitesimal-symprod}\@.  We let
$\gamma_1$ and $\gamma_2$ denote geodesics with initial conditions $X_1(x)$
and $X_2(x)$, respectively. Now define two new curves $\Upsilon_3^Z$ and
$\Upsilon_4^Z$ in $TM$ as follows:
\begin{align*}
\Upsilon_3^Z(t)=&\;\Phi^{X_2^V}_{-t}\tau_{\gamma_1}^{(0,t)}
\Phi^{X_2^V}_{t}\tau_{\gamma_1}^{(t,0)}\Phi^{X_1^V}_{-t}
\tau_{\gamma_2}^{(0,t)}\Phi^{X_1^V}_{t}\tau_{\gamma_2}^{(t,0)}(v_x),\\
\Upsilon_4^Z(t)=&\;\Phi^{X_2^V}_{-t}\overline{\tau}_{\gamma_1}^{(0,t)}
\Phi^{X_2^V}_{t}\overline{\tau}_{\gamma_1}^{(t,0)}\Phi^{X_1^V}_{-t}
\overline{\tau}_{\gamma_2}^{(0,t)}\Phi^{X_1^V}_{t}
\overline{\tau}_{\gamma_2}^{(t,0)}(v_x).
\end{align*}
With these constructions, we have the following result.
\begin{theorem}\label{the:infinitesimal-symprod-geodesics}
With the notation of the preceding paragraph, if\/
$\Upsilon\in\{\Upsilon_3^Z,\Upsilon_4^Z\}$\@, then\/
$\Upsilon'(0)=0$ and
\begin{equation*}
\frac{1}{2}\left.\frac{\rmd^2}{\rmd t^2}\right|_{t=0}\Upsilon(t)=
\langle X_1:X_2\rangle^V(v_x).
\end{equation*}
\begin{proof}
Corollary~\ref{cor:simple-infinitesimal} also applies to $\Upsilon\in
\{\Upsilon_3^Z,\Upsilon_4^Z\}$ by replacing $\eta_1$ and $\eta_2$ by
geodesics $\gamma_1$ and $\gamma_2$, respectively. From
Corollary~\ref{cor:simple-infinitesimal} it is easy to see that the first
derivative of $\Upsilon_3^Z$ at $t=0$ is zero. The second derivative at $t=0$
is
\begin{align*}
\left.\frac{\rmd^2 } {\rmd t^2}\right|_{t=0}\Upsilon_3^Z(t)=&\,
2\textup{vlft}\left(v_x,\left.\frac{\rmd}{\rmd t}\right|_{t=0}
(\tau^{(0,t)}_{\gamma_1}(X_2(\gamma_1(t))))+\left.\frac{\rmd}{\rmd t}
\right|_{t=0}(\tau^{(0,t)}_{\gamma_2}(X_1(\gamma_2(t))))\right)\\
=&\,2\textup{vlft}\left(v_x,\nabla_{\gamma_1'(0)}X_2(x)+
\nabla_{\gamma_2'(0)}X_1(x)\right),
\end{align*}
using~\eqref{eq:covartau}. Note that the last expression only depends on the
values and the first derivatives of the geodesics at zero. Since
$\gamma_1'(0)=X_1(x)$ and $\gamma_2'(0)=X_2(x)$, the theorem follows for
$\Upsilon_3^Z$.

The result can be proved for $\Upsilon_4^Z$ using
Lemma~\ref{lem:PT-SameGeod-DiffT} since the parallel transport in
$\Upsilon_4^Z$ is defined along geodesics.  If $\gamma$ is a geodesic and
$V\in T_{\gamma(t)}M$\@, then \eqref{eq:Diff-Paral-Trans} can be rewritten as
follows:
\begin{equation}\label{eq:Diff-Paral-Trans-Rewritten}
\overline{\tau}^{(0,t)}_\gamma(V)
=\tau^{(0,t)}_\gamma(V)+\frac{1}{2}\overline{\tau}_\gamma^{(0,t)}
\left(\int^0_t\overline{\tau}_\gamma^{(t,s)}
\left(T(\gamma'(s),\tau_\gamma^{(s,t)}(V))\right)\,
\rmd s\right),
\end{equation}
By Corollary~\ref{cor:simple-infinitesimal}
and~\eqref{eq:Diff-Paral-Trans-Rewritten} we have
\begin{multline*}
\Upsilon_4^Z(t)=\Upsilon_3^Z(t)+
t\left(\frac{1}{2}\overline{\tau}_{\gamma_2}^{(0,t)}\left(\int^0_t
\overline{\tau}_{\gamma_2}^{(t,s)}\left(
T(\gamma_2'(s),\tau_{\gamma_2}^{(s,t)}(X_1(\gamma_2(t))))\right)\,
\rmd s\right)\right.\\
\left.+\frac{1}{2}\overline{\tau}_{\gamma_1}^{(0,t)}
\left(\int^0_t\overline{\tau}_{\gamma_1}^{(t,s)}
\left(T(\gamma_1'(s),\tau_{\gamma_1}^{(s,t)}(X_2(\gamma_1(t))))\right)\,
\rmd s\right)\right).
\end{multline*}
As a result, with the abbreviation
\begin{multline*}
A(t)=\frac{1}{2}\overline{\tau}_{\gamma_2}^{(0,t)}\left(\int^0_t
\overline{\tau}_{\gamma_2}^{(t,s)}\left(T(\gamma_2'(s),
\tau_{\gamma_2}^{(s,t)}(X_1(\gamma_2(t))))\right)\,\rmd s\right)\\
+\frac{1}{2}\overline{\tau}_{\gamma_1}^{(0,t)}
\left(\int^0_t\overline{\tau}_{\gamma_1}^{(t,s)}
\left(T(\gamma_1'(s),\tau_{\gamma_1}^{(s,t)}(X_2(\gamma_1(t))))\right)\,
\rmd s\right),
\end{multline*}
we have
\begin{gather*}
\Upsilon_4^Z(0)= \Upsilon_3^Z(0), \qquad
\frac{\rmd }{\rmd t}\Upsilon_4^Z(t)=
\frac{\rmd }{\rmd t}\Upsilon_3^Z(t)+A(t)+t\frac{\rmd }{\rmd t}A(t),\\
\left.\frac{\rmd }{\rmd t}\right|_{t=0}\Upsilon_4^Z(t)=
\left.\frac{\rmd }{\rmd t}\right|_{t=0}\Upsilon_3^Z(t)+A(0)=
\left.\frac{\rmd }{\rmd t}\right|_{t=0}\Upsilon_1(t)=0,
\end{gather*}
since trivially $A(0)=0$\@.

Now
\begin{equation*}
\frac{\rmd^2}{\rmd t^2}\Upsilon_4^Z(t)=
\frac{\rmd^2 }{\rmd t^2}\Upsilon_3^Z(t)+2\frac{\rmd }{\rmd t}A(t)+
t\frac{\rmd^2 }{\rmd t^2}A(t).
\end{equation*}
At $t=0$,
\begin{equation*}
 \left.\frac{\rmd^2}{\rmd t^2}\right|_{t=0}\Upsilon_4^Z(t)=
\left.\frac{\rmd^2 }{\rmd t^2}\right|_{t=0}\Upsilon_3^Z(t)+
2\left.\frac{\rmd }{\rmd t}\right|_{t=0}A(t).
\end{equation*}
Note that
\begin{equation*}
\left.\frac{\rmd }{\rmd t}\right|_{t=0}A(t)=-T(X_2(x),X_1(x))-
T(X_1(x),X_2(x))=0.
\end{equation*}
Thus the result follows for $\Upsilon_4^Z$.
\end{proof}
\end{theorem}

\section{Characterization of geodesically invariant
distributions}\label{sec:geodesic-invariance}

In the preceding section we provided an interpretation of the symmetric
product that is similar to the composition of flows
formula~\eqref{eq:Lie-bracket} for the Lie bracket.  In this section we
provide an interpretation of the symmetric product rather like that which
Frobenius's Theorem provides for the Lie bracket.  The theorem we prove here
has already appeared in~\cite{ADL:96b,ADL:96a}\@.  However, we provide a
proof that is somewhat more elegant and also builds upon some independently
interesting constructions using distributions.

\subsection{Constructions using distributions}\label{subsec:distributions}

Let $M$ be a $n$-dimensional manifold with $\mathcal{D}$ a distribution on
$M$\@.  Distributions in this paper will always be smooth and of locally
constant rank.  Let $\tau_{\mathcal{D}}\colon TM/\mathcal{D}\to M$ be the
quotient vector bundle, and let $\pi_{\mathcal{D}}\colon TM\to
TM/\mathcal{D}$ be the canonical projection.  Note that the following diagram
commutes:
\begin{equation*}
\xymatrix{{TM}\ar[rr]^{\pi_{\mathcal{D}}}\ar[rd]_{\tau_M}&&
TM/\mathcal{D}\ar[ld]^{\tau_{\mathcal{D}}}\\&M&}
\end{equation*}

At $0_x\in TM/\mathcal{D}$\@, there exists the following natural splitting
\begin{equation}\label{eq:Split-0}
T_{0_x}TM/\mathcal{D}\simeq T_xM \oplus (T_xM/\mathcal{D}_x),
\end{equation}
cf.~\cite[Lemma~6.33]{FB/ADL:04}\@.  Hence we define the projection
$\textup{pr}_2\colon T_{0_x}TM/\mathcal{D} \rightarrow T_xM/\mathcal{D}_x$
onto the second component of the splitting in~\eqref{eq:Split-0}\@.  The
projection onto the first factor is simply $T_{0_x}\tau_{\mathcal{D}}$\@.

In the following result we give a characterization of vector fields tangent
to subbundles that will be useful for us.
\begin{proposition}\label{prop:Y-tangent-D-general}
Let\/ $\sigma\colon E\to M$ and\/ $\tau\colon F\to M$ be vector bundles and
let\/ $f\colon E\rightarrow F$ be a surjective vector bundle morphism over
the identity.  A vector field\/ $Y$ on\/ $E$ is tangent to\/ $\ker(f)$ if and
only if
\begin{equation*}
\textup{pr}_2(T_{e_x}f\circ Y(e_x))=0_x
\end{equation*}
for every\/ $e_x\in\ker(f)$\@.
\begin{proof}
For $x\in M$ the isomorphism of $T_{0_x}F$ with $T_xM\oplus F_x$ is given
explicitly by
\begin{equation*}
X_{0_x}\mapsto(T_{0_x}\tau(X_{0_x}),\textup{pr}_2(X_{0_x})).
\end{equation*}
Now note that, thinking of $\ker(f)=f^{-1}(Z(F))$ ($Z(F)$ is the zero section
of $F$ regarded as a submanifold of $F$) as a submanifold of $E$ we have, for
each $e_x\in\ker(f)$\@,
\begin{equation*}
T_{e_x}\ker(f)=\{X_{e_x}\in T_{e_x}E|\enspace
T_{e_x}f(X_{e_x})\in T_{0_x}Z(F)\}
\end{equation*}
(see~\cite[Theorem~3.5.12]{RA/JEM/TSR:88}).  Since
$T_{0_x}Z(F)=\textup{image}(T_{0_x}\tau)$ we have that $X_{e_x}\in
T_{e_x}\ker(f)$ if and only if
\begin{equation*}
\textup{pr}_2(T_{e_x}f(X_{e_x}))=0,
\end{equation*}
as desired.
\end{proof}
\end{proposition}

The following result is a particular case of
Proposition~\ref{prop:Y-tangent-D-general}\@, noting that
$\mathcal{D}=\ker(\pi_{\mathcal{D}})$\@.
\begin{corollary}\label{cor:Y-tangent-D-final}
Let\/ $\mathcal{D}\subset TM$ be a distribution.  A vector field\/ $Y$ on\/
$TM$ is tangent to\/ $\mathcal{D}$ if and only if
\begin{equation*}
\textup{pr}_2((T_{v_x}\pi_{\mathcal{D}}\circ Y)(v_x))=0_x
\end{equation*}
for every\/ $v_x\in\mathcal{D}$\@.
\end{corollary}

A corollary to this corollary, and one that will be useful for us, is the
following.
\begin{corollary}\label{cor:vlft-dist}
A vector field\/ $X\in\Gamma^\infty(TM)$ takes values in a distribution\/
$\mathcal{D}$ on\/ $M$ if and only if\/ $X^V$ is tangent to\/
$\mathcal{D}$\@.
\begin{proof}
First suppose that $X$ is $\mathcal{D}$-valued.  If $v_x\in\mathcal{D}$ then,
for any $t\in\mathbb{R}$\@,
\begin{equation*}
v_x+tX(x)\in\mathcal{D}\enspace\implies\enspace
\pi_{\mathcal{D}}(v_x+tX(x))=0_x
\enspace\implies\enspace T_{v_x}\pi_{\mathcal{D}}(X^V(v_x))=0,
\end{equation*}
giving $X^V$ tangent to $\mathcal{D}$ by the previous corollary.

Conversely, suppose that $X^V$ is tangent to $\mathcal{D}$\@.  Then, by the
previous corollary,
\begin{equation*}
\textup{pr}_2(T_{v_x}\pi_{\mathcal{D}}(X^V(v_x)))=0_x
\end{equation*}
for every $v_x\in\mathcal{D}$\@.  Since $X^V$ is vertical and since
$\pi_{\mathcal{D}}$ is a vector bundle mapping,
$T_{v_x}\pi_{\mathcal{D}}\circ X^V$ is vertical.  Thus
\begin{equation*}
T_{0_x}\tau_{\mathcal{D}}(T_{v_x}\pi_{\mathcal{D}}\circ X^V)(v_x)=0_x.
\end{equation*}
This implies that both components of $T_{v_x}\pi_{\mathcal{D}}\circ X^V(v_x)$
are zero in the decomposition $T_{0_x}TM/\mathcal{D}\simeq
T_xM\oplus(T_xM/\mathcal{D}_x)$\@.  Thus, reversing the calculations from the
first part of the proof, we conclude that $X(x)\in\mathcal{D}$\@.
\end{proof}
\end{corollary}

\subsection{The geodesic invariance theorem}

Let us define the objects of interest.
\begin{definition}
A distribution $\mathcal{D}$ on $M$ is \textbf{geodesically invariant} under
an affine connection $\nabla$ on $M$ if, for every geodesic $\gamma\colon
I\to M$ for which $\gamma'(t_0)\in\mathcal{D}_{\gamma(t_0)}$ for some $t_0\in
I$\@, it holds that $\gamma'(t)\in\mathcal{D}_{\gamma(t)}$ for every $t\in
I$\@.
\end{definition}

Before we prove the main result in this section, we need to prove a few
technical lemmata.  The first relates the horizontal lift of a vector field
to the complete lift of the same vector field.  A special case of this
formula is given by~\cite{AB:07,AB:10}\@, with an intrinsic proof using frame
bundles in~\cite{AB:10}\@.
\begin{lemma}\label{lem:XC-XH}
If\/ $X\in\Gamma^\infty(TM)$ then
\begin{equation*}
X^C(v_x)=X^H(v_x)+_1\textup{vlft}(v_x,\nabla_{v_x}X(x)+T(X(x),v_x)),
\end{equation*}
for every\/ $v_x\in TM$\@.
\begin{proof}
A direct proof in coordinates is, of course, elementary.  However, we shall
provide an intrinsic proof to keep in the spirit of our intrinsic proof of
Theorem~\ref{the:Geod-Inv} below.

Let $v_x\in TM$ and let $Y\in\Gamma^\infty(TM)$ be such that $Y(x)=v_x$\@.
Note that
\begin{equation*}
\left.\frac{\rmd}{\rmd s}\right|_{s=0}\Phi^X_t\Phi^Y_s(x)=
T_x\Phi^X_t(Y(x)).
\end{equation*}
Also compute
\begin{align*}
\left.\frac{\rmd}{\rmd s}\right|_{s=0}\Phi^X_t\Phi^Y_s(x)=&\;
\left.\frac{\rmd}{\rmd s}\right|_{s=0}\Phi^Y_s\Phi^X_t\Phi^X_{-t}
\Phi^Y_{-s}\Phi^X_t\Phi^Y_s(x)\\
=&\;Y(\Phi^X_t(x))+T_x\Phi^X_t\left(\left.\frac{\rmd}{\rmd s}\right|_{s=0}
\Phi^X_{-t}\Phi^Y_{-s}\Phi^X_t\Phi^Y_s(\Phi^X_t(x))\right).
\end{align*}
Note that
\begin{align*}
\textup{BCH}_1(sY,tX,-sY,-tX)=&\;0,\\
\textup{BCH}_2(sY,tX,-sY,-tX)=&\;st[Y,X].
\end{align*}
Therefore, using~\eqref{eq:BCH12}\@,
\begin{align*}
\left.\frac{\rmd}{\rmd s}\right|_{s=0}
\Phi^X_{-t}\Phi^Y_{-s}\Phi^X_t\Phi^Y_s(\Phi^X_t(x))=&\;
\left.\frac{\rmd}{\rmd s}\right|_{s=0}\Phi^{st[Y,X]}_1(\Phi^X_t(x))\\
=&\;\left.\frac{\rmd}{\rmd
s}\right|_{s=0}\Phi^{t[Y,X]}_s(\Phi^X_t(x))= t[Y,X](\Phi^X_t(x)).
\end{align*}
Putting the above calculations together gives
\begin{equation*}
T_x\Phi^X_t(Y(x))=Y(\Phi^X_t(x))-t[X,Y](\Phi^X_t(x)).
\end{equation*}
Thus, recalling~\eqref{eq:XC-flow}\@,
\begin{equation*}
\Phi^{-X^H}_t\Phi^{X^C}_t(Y(x))=
\tau^{(t,0)}_{\gamma_-}(Y(\Phi^X_t(x))-t[X,Y](\Phi^X_t(x))),
\end{equation*}
where $\gamma_-$ is the integral curve of $-X$ through $\Phi^X_t(x)$\@, where
we have used~\eqref{eq:hlft-flow}\@.  If $\gamma$ is the integral curve of
$X$ through $x$ note that $\tau^{(t,0)}_{\gamma_-}=\tau^{(0,t)}_\gamma$\@.
Now we compute
\begin{align*}
\left.\frac{\rmd}{\rmd
t}\right|_{t=0}\Phi^{-X^H}_t\Phi^{X^C}_t(Y(x))=&\;
\left.\frac{\rmd}{\rmd
t}\right|_{t=0}\tau^{(0,t)}_\gamma(Y(\Phi^X_t(x))-
t[X,Y](\Phi^X_t(x)))\\
=&\;\nabla_XY(x)-[X,Y](x)=\nabla_YX(x)+T(X(x),Y(x)),
\end{align*}
using~\eqref{eq:covartau}\@.  Note that since $X^C$ and $X^H$ are both vector
fields over $X$\@, it follows that
\begin{equation*}
t\mapsto\tau^{(0,t)}_\gamma(Y(\Phi^X_t(x)))
\end{equation*}
is a curve in $T_xM$\@.  Thus the derivative of this curve at $t=0$ is in
$V_{Y(x)}TM$\@.  Thus we have shown that
\begin{equation}\label{eq:XHXC}
\left.\frac{\rmd}{\rmd t}\right|_{t=0}\Phi^{-X^H}_t\Phi^{X^C}_t(v_x)=
\textup{vlft}(v_x,\nabla_{v_x}X(x)+T(X(x),v_x)).
\end{equation}

Finally, by the BCH formula, we have
\begin{equation*}
\Phi^{-X^H}_t\Phi^{X^C}_t(v_x)=
\Phi^{\phi(\textup{BCH}_1(-tX^H,tX^C))}_1(v_x)+O(|t|^2)=
\Phi^{X^C-X^H}_t(v_x)+O(|t|^2).
\end{equation*}
Differentiating with respect to $t$ and evaluating at $t=0$\@,
using~\eqref{eq:XHXC}\@, gives the result.
\end{proof}
\end{lemma}

Another useful lemma is the following.
\begin{lemma}\label{lem:XVZYV}
If\/ $Z$ is the geodesic spray for an affine connection and if\/
$X,Y\in\Gamma^\infty(TM)$\@, then
\begin{equation*}
[X^V,[Z,Y^V]]=\langle X\colon Y\rangle^V.
\end{equation*}
\begin{proof}
Again, a proof in coordinates is easy, but we give an intrinsic proof.

We use the following formula for the Lie
bracket~\cite[Theorem~4.2.19]{RA/JEM/TSR:88}\@:
\begin{equation*}
[U,V](v_x)=\left.\frac{\rmd}{\rmd t}\right|_{t=0}(\Phi^U_t)^*V(v_x),
\end{equation*}
for vector fields $U$ and $V$ on $TM$\@.  Note that
$\Phi^{Y^V}_t=\textup{Id}_{TM}+tY\circ\tau_M$\@,
using~\eqref{eq:vlft-flow}\@, and so
\begin{equation*}
Z\circ\Phi^{Y^V}_t(v_x)=\textup{hlft}(v_x+tY(x),v_x+tY(x)),
\end{equation*}
using~\eqref{eq:Zdefn}\@.  Now note that for $U\in\Gamma^\infty(TM)$\@,
\begin{equation*}
T\tau_M=T\tau_M\circ TU\circ T\tau_M,
\end{equation*}
using the fact that $\tau_M\circ U=\textup{Id}_M$\@.  It follows that, if
$W\in TTM$\@, the expression $W+_2TU\circ T\tau_M(W)$ makes sense.  Thus we
have
\begin{equation*}
T\Phi^{Y^V}_{-t}(W)=W-_2T(tY)\circ T\tau_M(W).
\end{equation*}
Thus
\begin{align*}
(\Phi^{Y^V}_t)^*Z(v_x)=&\;T\Phi^{Y^V}_{-t}\circ Z\circ\Phi^{Y^V}_t(v_x)\\
=&\;\textup{hlft}(v_x+tY(x),v_x+tY(x))-_2T(tY)(v_x+tY(x)).
\end{align*}
We need to differentiate this expression with respect to $t$\@.  To do this,
let us define $\Upsilon\colon\mathbb{R}^2\to TTM$ by
\begin{equation*}
\Upsilon(s,t)=\textup{hlft}(v_x+sY(x),v_x+tY(x))-_2T(sY)(v_x+tY(x))
\end{equation*}
and $\iota\colon\mathbb{R}\to\mathbb{R}^2$ by $\iota(t)=(t,t)$\@.  Note that
\begin{equation*}
(\Phi^{Y^V}_t)^*Z(v_x)=\Upsilon\circ\iota(t)
\end{equation*}
and so
\begin{equation*}
\left.\frac{\rmd}{\rmd t}\right|_{t=0}(\Phi^{Y^V}_t)^*Z(v_x)=
T_1\Upsilon(1,1)+_1T_2\Upsilon(1,1),
\end{equation*}
where $T_1\Upsilon$ and $T_2\Upsilon$ denote the partial derivatives of
$\Upsilon$\@,~cf.~\cite[Proposition~3.3.13]{RA/JEM/TSR:88}\@.  Thus we have
\begin{equation*}
T_1\Upsilon(1,1)=\left.\frac{\rmd}{\rmd s}\right|_{s=0}\Upsilon(s,0),\quad
T_2\Upsilon(1,1)=\left.\frac{\rmd}{\rmd t}\right|_{t=0}\Upsilon(0,t).
\end{equation*}
The second of these expressions is readily calculated:
\begin{align*}
\left.\frac{\rmd}{\rmd t}\right|_{t=0}\Upsilon(0,t)
=&\;\left.\frac{\rmd}{\rmd t}\right|_{t=0}\textup{hlft}(v_x,v_x+tY(x))\\
=&\;\left.\frac{\rmd}{\rmd t}\right|_{t=0}\bigl(\textup{hlft}(v_x,v_x)
+_1t\cdot_1\textup{hlft}(v_x,Y(x))\bigr)\\
=&\;\textup{hlft}(v_x,Y(x)),
\end{align*}
using~\eqref{eq:Hdiagrams}\@.  For the first, note that
\begin{equation*}
\tau_{TM}(\textup{hlft}(v_x+sY(x),v_x)-_2T(sY)(v_x))=v_x,\quad
T\tau_M(\textup{hlft}(v_x+sY(x),v_x)-_2T(sY)(v_x))=v_x.
\end{equation*}
Hence the tangent vector to the curve $0\mapsto \Upsilon(s,0)$ at every time
$s$ lies in the vertical subspace $V_{\Upsilon(s,0)}(TTM)\simeq
T_{v_x}TM\simeq T_xM \oplus T_xM$\@.  In fact, to be more precise, it lies in
the second copy of $T_xM$\@.  In other words
\begin{equation*}
\frac{\rmd}{\rmd s}\Upsilon(s,0)\in
\{\textup{vlft}(\Upsilon(s,0),\textup{vlft}(0_x,w_x))|\enspace
w_x\in T_xM\}\simeq 0_x\oplus T_xM\subset V_{\Upsilon(s,0)}(TTM),
\end{equation*}
where the first $\textup{vlft}$ lifts from $TTM$ to $VTTM$ and the
second lifts from $TM$ to $VTM$\@.  Thus,
\begin{equation*}
\left.\frac{\rmd}{\rmd s}\right|_{s=0}
\bigl(\textup{hlft}(v_x+sY(x),v_x)-_2T(sY)(v_x)\bigr)=
\textup{vlft}(\textup{hlft}(v_x,v_x),\textup{hlft}(Y(x),v_x)-_2TY(v_x)).
\end{equation*}
Then, using the definition of $\textup{hlft}$ from~\eqref{eq:Hdef}
and Lemma~\ref{lem:intertwine} we obtain
\begin{align*}
\textup{hlft}(Y(x),v_x)-_2TY(v_x)=&\;(TY(v_x)-_1
\textup{vlft}(Y(x),\nabla_vY(x)))-_2TY(v_x)\\
=&\;(TY(v_x)+_1\textup{vlft}(Y(x),-\nabla_vY(x)))-_2(TY(v_x)+_1\xi_1(Y(x)))\\
=&\;(TY(v_x)-_2TY(v_x))+_1(\textup{vlft}(Y(x),-\nabla_vY(x))-_2\xi_1(Y(x)))\\
=&\;\textup{vlft}(0_x,-\nabla_vY(x)),
\end{align*}
where $\xi_1\colon TM\rightarrow TTM$ is the zero section relative to the
primary vector bundle structure.  According to the identification mentioned
above, we have
\begin{multline*}
\left.\frac{\rmd}{\rmd s}\right|_{s=0}
\bigl(\textup{hlft}(v_x+sY(x),v_x)-_2T(sY)(v_x)\bigr)\\=
\textup{vlft}(\textup{hlft}(v_x,v_x),\textup{vlft}(0_x,-\nabla_vY(x)))\simeq
-\textup{vlft}(v_x,\nabla_vY(x)).
\end{multline*}

Putting the above together gives
\begin{equation*}
[Y^V,Z](v_x)=\textup{hlft}(v_x,Y(x))-_1\textup{vlft}(v_x,\nabla_{v_x}Y).
\end{equation*}

In like manner we compute
\begin{multline*}
(\Phi^{X^V}_t)^*[Y^V,Z](v_x)=\\\bigl(\textup{hlft}(v_x+tX(x),Y(x))
-_1\textup{vlft}(v_x+tX(x),\nabla_{v_x+tX(x)}Y)\bigr)-_2T(tX)(Y(x)).
\end{multline*}
We differentiate this expression as above, in this case defining
\begin{equation*}
\Upsilon(s,t)=\bigl(\textup{hlft}(v_x+sX(x),Y(x))
-_1\textup{vlft}(v_x+sX(x),\nabla_{v_x+tX(x)}Y)\bigr)-_2T(sX)(Y(x)).
\end{equation*}
The two expressions we need to differentiate are then
\begin{equation*}
\Upsilon(s,0)=\bigl(\textup{hlft}(v_x+sX(x),Y(x))
-_1\textup{vlft}(v_x+sX(x),\nabla_{v_x}Y)\bigr)-_2T(sX)(Y(x))
\end{equation*}
and
\begin{equation*}
\Upsilon(0,t)=\bigl(\textup{hlft}(v_x,Y(x))
-_1\textup{vlft}(v_x,\nabla_{v_x+tX(x)}Y)\bigr).
\end{equation*}
The second of these is easily differentiated:
\begin{equation*}
\left.\frac{\rmd}{\rmd t}\right|_{t=0}\Upsilon(0,t)=
-\textup{vlft}(v_x,\nabla_XY(x)),
\end{equation*}
using~\eqref{eq:Vdiagrams}\@.  For the first, we first note that
\begin{multline*}
\left.\frac{\rmd}{\rmd s}\right|_{s=0}\bigl(\textup{hlft}(v_x+sX(x),Y(x))
-_1\textup{vlft}(v_x+sX(x),\nabla_{v_x}Y)\bigr)-_2T(sX)(Y(x))\\
=\left.\frac{\rmd}{\rmd s}\right|_{s=0}
\textup{hlft}(v_x+sX(x),Y(x))-_2T(sX)(Y(x))
\end{multline*}
since the vertical component of the second term is independent of $s$\@.  Now
we can proceed as above to compute
\begin{multline*}
\left.\frac{\rmd}{\rmd s}\right|_{s=0}\Upsilon(s,0)=
\left.\frac{\rmd}{\rmd s}\right|_{s=0}
\bigl(\textup{hlft}(v_x+sX(x),Y(x))-_2T(sX)(Y(x))\bigr)\\
=\textup{vlft}(\textup{hlft}(v_x,Y(x)),\textup{hlft}(X(x),Y(x))-_2TX(Y(x))).
\end{multline*}
Then, as above, using the definition of $\textup{hlft}$
from~\eqref{eq:Hdef} and Lemma~\ref{lem:intertwine} we obtain
\begin{align*}
\textup{hlft}(X(x),&Y(x))-_2TX(Y(x))\\
=&\;(TX(Y(x))-_1\textup{vlft}(X(x),\nabla_YX(x)))-_2
(TX(Y(x))+_1\xi_1(X(x)))\\
=&\;\textup{vlft}(X(x),-\nabla_YX(x))-_2\xi_1(X(x))\\
=&\;\textup{vlft}(0_x,-\nabla_YX(x)).
\end{align*}
According to the identification mentioned above, we have
\begin{multline*}
\left.\frac{\rmd}{\rmd s}\right|_{s=0}
\bigl(\textup{hlft}(v_x+sY(x),v_x)-_2T(sY)(v_x)\bigr)\\=
\textup{vlft}(\textup{hlft}(v_x,Y(x)),\textup{vlft}(0_x,-\nabla_YX(x)))\simeq
-\textup{vlft}(v_x,\nabla_XY(x)).
\end{multline*}

Putting the preceding calculations together and appropriately
identifying vertical tangent vectors gives
\begin{equation*}
[X^V,[Y^V,Z]](v_x)=-\textup{vlft}(v_x,\nabla_XY(x))-_1
\textup{vlft}(v_x,\nabla_YX(x))=-\textup{vlft}(v_x,\langle X\colon
Y\rangle(x)),
\end{equation*}
which is the result.
\end{proof}
\end{lemma}

For $X\in\Gamma^\infty(TM)$ let us denote
\begin{equation*}
\mathcal{D}_X=\{\alpha X(x)|\enspace x\in M,\ \alpha\in\mathbb{R}\}.
\end{equation*}
With this notation, the last technical lemma upon which we shall draw is the
following.
\begin{lemma}\label{lem:XH-restricted}
A distribution\/ $\mathcal{D}$ is geodesically invariant if and only if, for
each\/ $X\in\Gamma^\infty(\mathcal{D})$ and for each\/
$v_x\in\mathcal{D}_X$\@, $X^H(v_x)\in T_{v_x}\mathcal{D}$\@.
\begin{proof}
First suppose that $X^H(v_x)\in T_{v_x}\mathcal{D}$ for every
$X\in\Gamma^\infty(\mathcal{D})$ and every $v_x\in\mathcal{D}_X$\@.  Let
$v_x\in\mathcal{D}$ and let $X\in\Gamma^\infty(\mathcal{D})$ be such that
$X(x)=v_x$\@.  (This is possible as follows.  Since $\mathcal{D}$ is smooth
and constant rank, there exists linearly independent smooth local generators
$X_1,\dots,X_k$ for $\mathcal{D}$ about $x$\@.  Write
\begin{equation*}
v_x=\alpha_1X_1(x)+\dots+\alpha_kX_k(x)
\end{equation*}
and let $f_1,\dots,f_k\colon M\to\mathbb{R}$ be such that $f_j(x)=\alpha_j$
and such that $f_1,\dots,f_k$ vanish outside a sufficiently small
neighbourhood of $x$\@.  Then take
\begin{equation*}
X=f_1X_1+\dots+f_kX_k.)
\end{equation*}
By hypothesis, $X^H(X(x))\in T_{v_x}\mathcal{D}$\@.  By~\eqref{eq:Zdefn} and
the definition of $X^H$ it follows that $Z(v_x)\in T_{v_x}\mathcal{D}$\@.  As
$v_x\in\mathcal{D}$ is arbitrary, it follows that $Z$ is tangent to
$\mathcal{D}$\@, meaning that $\mathcal{D}$ is geodesically invariant.

Conversely, suppose that $Z(v_x)\in T_{v_x}\mathcal{D}$ for every
$v_x\in\mathcal{D}$\@.  Let $X\in\Gamma^\infty(\mathcal{D})$ and let
$v_x\in\mathcal{D}_X$\@.  Thus $v_x=\alpha X(x)$ for some
$\alpha\in\mathbb{R}$\@.  We then have
\begin{equation*}
T_{v_x}\mathcal{D}\ni Z(v_x)=\textup{hlft}(v_x,v_x)=
\alpha\textup{hlft}(v_x,X(x))=\alpha X^H(v_x),
\end{equation*}
using~\eqref{eq:Zdefn}\@.  We then consider two cases.  First of all, suppose
that $X(x)=0_x$\@.  Then $v_x=0_x$ and so $Z(v_x)=X^H(v_x)=0_{v_x}$ and we
trivially have $X^H(v_x)\in T_{v_x}\mathcal{D}$\@.  If $X(x)\not=0_x$ then
our computation just preceding gives $X^H(v_x)=\alpha^{-1}Z(v_x)\in
T_{v_x}\mathcal{D}$\@.
\end{proof}
\end{lemma}

We can now state the main result in this section.  While this result is
known~\cite{ADL:96b,ADL:96a}\@, we provide here a self-contained intrinsic
proof using the tools developed in the paper.
\begin{theorem}[\cite{ADL:96b,ADL:96a}]\label{the:Geod-Inv}
Let\/ $\mathcal{D}$ be a distribution on a manifold\/ $M$ with an affine
connection\/ $\nabla$\@.  The following are equivalent:
\begin{compactenum}
\item \label{pl:geoinv1} $\mathcal{D}$ is geodesically invariant;
\item \label{pl:geoinv2} $\langle X\colon Y\rangle\in\Gamma^\infty(\mathcal{D})$ for every\/
$X,Y\in\Gamma^\infty(\mathcal{D})$;
\item \label{pl:geoinv3} $\nabla_XX\in\Gamma^\infty(\mathcal{D})$ for every\/
$X\in\Gamma^\infty(\mathcal{D})$\@.
\end{compactenum}
\begin{proof}
\eqref{pl:geoinv1}$\implies$\eqref{pl:geoinv2} The proof of this
in~\cite{ADL:96b,ADL:96a} makes use of the formula from Lemma~\ref{lem:XVZYV}
which was only derived there in coordinates.  We reproduce this proof here,
but now it is a self-contained intrinsic proof since we have an intrinsic
proof of Lemma~\ref{lem:XVZYV}\@.  We also provide a second proof using our
composition formula from Theorem~\ref{the:infinitesimal-symprod-geodesics}
for the symmetric product.

\textit{First proof:} Let $X,Y\in\Gamma^\infty(\mathcal{D})$\@.  It is clear
that since $\mathcal{D}$ is geodesically invariant, $Z$ is tangent to
$\mathcal{D}$\@.  Moreover, by Corollary~\ref{cor:vlft-dist}\@, $X^V$ and
$Y^V$ are tangent to $\mathcal{D}$\@.  By the
formula~\eqref{eq:Lie-bracket}\@, it follows that all Lie brackets involving
$Z$\@, $X^V$\@, and $Y^V$ are also tangent to $\mathcal{D}$\@.  In
particular, $[X^V,[Z,Y^V]]$ is tangent to $\mathcal{D}$ and so, by
Corollary~\ref{cor:vlft-dist} and Lemma~\ref{lem:XVZYV}\@, $\langle X\colon
Y\rangle$ is tangent to $\mathcal{D}$\@.

\textit{Second proof:} By
Theorem~\ref{the:infinitesimal-symprod-geodesics} we know that
\begin{equation*}
\frac{1}{2}\left.\frac{{\rm d}^2}{{\rm d}t^2}\right|_{t=0}
\Upsilon^Z_3(t)=\langle X\colon Y \rangle^V(v_x).
\end{equation*}
In particular this is true for every $X,Y\in\Gamma^\infty(\mathcal{D})$\@.
According to Corollary~\ref{cor:Y-tangent-D-final} we only have to prove that
\begin{equation*}
{\rm pr}_2\left(T_{v_x}\pi_{\mathcal{D}} \circ \left.\frac{{\rm
d}^2}{{\rm d}t^2}\right|_{t=0} \Upsilon^Z_3(t)\right)=0_x
\end{equation*}
for every $v_x\in \mathcal{D}$\@.  First assume that $X=Y$\@.  By adapting
conveniently Corollary~\ref{cor:simple-infinitesimal} to
$\Upsilon_3^Z$ we have
\begin{equation*}
\frac{1}{2} \, \left.\frac{{\rm d}^2}{{\rm d}t^2}\right|_{t=0}
\Upsilon^Z_3(t)= \left.\frac{\rm d}{{\rm d}t}\right|_{t=0}
\tau^{(0,t)}_{\gamma_X}(X(\gamma_X(t))),
\end{equation*}
where $\gamma_X$ is a geodesic such that $\gamma_X(0)=x$ and
$\gamma_X'(0)=X(x)$\@.  By the Leibniz rule,
\begin{equation*}
\left.\frac{\rm d}{{\rm d}t}\right|_{t=0}
\tau^{(0,t)}_{\gamma_X}(X(\gamma_X(t)))=
\left(\left.\frac{\rmd}{\rmd t}\right|_{t=0}
\tau^{(0,t)}_{\gamma_X}\right)(X(x))=-Z(X(x)).
\end{equation*}
As $\mathcal{D}$ is geodesically invariant by hypothesis, $Z$ is
tangent to $\mathcal{D}$\@. Using
Corollary~\ref{cor:Y-tangent-D-final} and the polarization identity
for the symmetric product, the result follows.

\eqref{pl:geoinv2}$\implies$\eqref{pl:geoinv3} This follows from the
definition of the symmetric product.

\eqref{pl:geoinv3}$\implies$\eqref{pl:geoinv1} Let
$X\in\Gamma^\infty(\mathcal{D})$ and let $v_x\in\mathcal{D}_X$\@.  Since
$X\in\Gamma^\infty(\mathcal{D})$\@, $\pi_{\mathcal{D}}\circ X(y)=0_y$ for
every $y\in M$\@.  Thus
\begin{equation*}
T_{X(x)}\pi_{\mathcal{D}}\circ T_xX(u_x)=u_x\oplus 0_x,
\end{equation*}
using the identification $T_{0_x}(TM/\mathcal{D})\simeq T_xM\oplus
T_xM/\mathcal{D}_x$\@.  This gives, in particular,
\begin{equation*}0_x=
\textup{pr}_2\circ T_{X(x)}\pi_{\mathcal{D}}\circ T_xX(X(x))=
\textup{pr}_2\circ T_{X(x)}\pi_{\mathcal{D}}(X^H(X(x))+(\nabla_XX)^V(X(x)))
\end{equation*}
using Lemma~\ref{lem:XC-XH}\@.  By hypothesis and by
Corollary~\ref{cor:vlft-dist}\@, $(\nabla_XX)^V$ is tangent to
$\mathcal{D}$\@.  Therefore, by Corollary~\ref{cor:Y-tangent-D-final}\@,
\begin{equation*}
\textup{pr}_2\circ
T_{X(x)}\pi_{\mathcal{D}}((\nabla_XX)^V(X(x)))=0_x.
\end{equation*}
Another appeal to Corollary~\ref{cor:Y-tangent-D-final} then allows
us to conclude that $X^H(X(x))\in T_{X(x)}\mathcal{D}$\@.  Linearity
of horizontal lift implies that $X^H(v_x)\in T_{v_x}\mathcal{D}$ for
all $v_x \in\mathcal{D}_X$\@, and the theorem follows from
Lemma~\ref{lem:XH-restricted}\@.
\end{proof}
\end{theorem}

\section*{Acknowledgements}
The authors are grateful to Professor Miguel C. Mu\~noz-Lecanda for
useful comments. MBL has been partially supported by MICINN (Spain)
Grants MTM2008-00689 and MTM2009-08166; 2009SGR1338 of the Catalan
government, IRSES project GEOMECH (246981) within the 7th European
Community Framework Program, by Beatriu de Pin\'os fellowship from
Comissionat per a Universitats i Recerca del Departament
d'Innovaci\'o, Universitats i Empresa of Generalitat de Catalunya
and by Juan de la Cierva fellowship from MICINN.
\bibliographystyle{abbrv}
\bibliography{refs}

\begin{thebibliography}{10}

\bibitem{RA/JEM/TSR:88}
R.~Abraham, J.~E. Marsden, and T.~S. Ratiu.
\newblock {\em Manifolds, Tensor Analysis, and Applications}.
\newblock Number~75 in Applied Mathematical Sciences. Springer-Verlag, 2
  edition, 1988.

\bibitem{MBL/MS:10}
M.~Barbero-Li{\~n}{\'a}n and M.~Sigalotti.
\newblock High-order sufficient conditions for configuration tracking of affine
  connection control systems.
\newblock {\em Systems {\&} Control Letters}, 59(8):491--503, 2010.

\bibitem{AB:07}
A.~Bhand.
\newblock {\em Geodesic Reduction via Frame Bundle Geometry}.
\newblock PhD thesis, Queen's University, Kingston, Department of Mathematics
  \& Statistics, Kingston, ON~K7L~3N6, Canada, June 2007.

\bibitem{AB:10}
A.~Bhand.
\newblock Geodesic reduction via frame bundle geometry.
\newblock {\em Symmetry, Integrability and Geometry: Methods and Applications},
  6(020):17 pages, 2010.

\bibitem{FB/NEL/ADL:97}
F.~Bullo, N.~E. Leonard, and A.~D. Lewis.
\newblock Controllability and motion algorithms for underactuated {L}agrangian
  systems on {L}ie groups.
\newblock {\em Institute of Electrical and Electronics Engineers. Transactions
  on Automatic Control}, 45(8):1437--1454, 2000.

\bibitem{FB/ADL:04}
F.~Bullo and A.~D. Lewis.
\newblock {\em Geometric Control of Mechanical Systems: Modeling, Analysis, and
  Design for Simple Mechanical Systems}.
\newblock Number~49 in Texts in Applied Mathematics. Springer-Verlag, New
  York\textendash{}Heidelberg\textendash{}Berlin, 2004.

\bibitem{MC:00}
M.~Crampin.
\newblock Connections of {B}erwald type.
\newblock {\em Universitatis Debreceniensis. Institutum Mathematicum.
  Publicationes Mathematicae}, 57(3-4):455--473, 2000.

\bibitem{PEC:81}
P.~E. Crouch.
\newblock Geometric structures in systems theory.
\newblock {\em Institution of Electrical Engineers. Proceedings. D. Control
  Theory and Applications}, 128(5):242--252, 1981.

\bibitem{RJF/HTL:99}
R.~J. Fisher and H.~T. Laquer.
\newblock Second order tangent vectors in {R}iemannian geometry.
\newblock {\em Journal of the Korean Mathematical Society}, 36(5):959--1008,
  1999.

\bibitem{SK/KN:63a}
S.~Kobayashi and K.~Nomizu.
\newblock {\em Foundations of Differential Geometry, {\normalfont Volume I}}.
\newblock Number~15 in Interscience Tracts in Pure and Applied Mathematics.
  Interscience Publishers, New York, 1963.

\bibitem{MK/JEM:10}
M.~Kobilarov and J.~E. Marsden.
\newblock Discrete geometric optimal control on {L}ie groups.
\newblock {\em Institute of Electrical and Electronics Engineers. Transactions
  on Robotics}, 2011.
\newblock To appear.

\bibitem{IK/PWM/JS:93}
I.~Kol{\'a\v{r}}, P.~W. Michor, and J.~Slov{\'a}k.
\newblock {\em Natural Operations in Differential Geometry}.
\newblock Springer-Verlag, New York\textendash{}Heidelberg\textendash{}Berlin,
  1993.

\bibitem{ADL:96b}
A.~D. Lewis.
\newblock A symmetric product for vector fields and its geometric meaning.
\newblock Technical Memorandum CIT--CDS 93--003, California Institute of
  Technology, Pasadena, CA 91125, 1996.

\bibitem{ADL:96a}
A.~D. Lewis.
\newblock Affine connections and distributions with applications to
  nonholonomic mechanics.
\newblock {\em Reports on Mathematical Physics}, 42(1/2):135--164, 1998.

\bibitem{ADL/RMM:95c}
A.~D. Lewis and R.~M. Murray.
\newblock Controllability of simple mechanical control systems.
\newblock {\em SIAM Journal on Control and Optimization}, 35(3):766--790, 1997.

\bibitem{KCHM:05}
K.~C.~H. Mackenzie.
\newblock {\em The General Theory of {L}ie Groupoids and {L}ie Algebroids}.
\newblock Number 213 in London Mathematical Society Lecture Note Series.
  Cambridge University Press, New York/Port Chester/Melbourne/Sydney, 2005.

\bibitem{JP:74}
J.~Pradines.
\newblock Repr\'esentation des jets non holonomes par des morphismes vectoriels
  doubles soud\'es.
\newblock {\em Comptes Rendus de l'Acad{\'e}mie des Sciences. S{\'e}rie I.
  Math{\'e}matique}, 278:1523--1526, 1974.

\bibitem{JPS:92}
J.-P. Serre.
\newblock {\em {L}ie Algebras and {L}ie Groups}.
\newblock Number 1500 in Lecture Notes in Mathematics. Springer-Verlag, New
  York\textendash{}Heidelberg\textendash{}Berlin, 1992.

\bibitem{RSS:87}
R.~S. Strichartz.
\newblock The
  {C}ampbell\textendash{}{B}aker\textendash{}{H}ausdorff\textendash{}{D}ynkin
  formula and solutions of differential equations.
\newblock {\em Journal of Functional Analysis}, 72:320--345, 1987.

\end{thebibliography}
\end{document}